\documentclass[%
 aip,
 amsmath,amssymb
 reprint,%
]{revtex4-1}
\usepackage{graphicx}
\usepackage{dcolumn}
\usepackage{bm}

\usepackage[dvipsnames]{xcolor}
\usepackage{soul}
\usepackage[utf8]{inputenc}
\usepackage[T1]{fontenc}
\usepackage{mathptmx}
\usepackage{etoolbox}


\begin{document}


\title{Bioconvection in a phototactic algae suspension with oblique irradiation and forward anisotropic scattering }
\author{Sandeep Kumar}
\email{sandeepkumar1.00123@gmail.com}
\author{Preeti Sharma}

\affiliation{
Department of Mathematics, PDPM  Indian Institute of Information Technology Design and Manufacturing, Jabalpur 482005, India
}%



\begin{abstract}
In this study, we analyze the bioconvection in a suspension of phototactic algae that exhibits anisotropic scattering. The top layer of the suspension is illuminated by oblique collimated irradiation. During the study, the bottom boundary is considered as rigid whereas the top boundary is considered stress-free. In order to solve the eigenvalue problem, the Newton-Raphson-Kantorovich finite difference method of order four is used. Linear analysis of the basic state is performed using neutral curves. The results demonstrate a change in the most unstable mode from an overstable to a stationary state or vice versa for particular parameters in response to a variation in the incidence angle. The position of the maximum basic concentration shifts toward the top of the suspension as the incidence angle is increased. In most cases, the system becomes more unstable with an increment in the incidence angle.

\end{abstract}

\maketitle

\section{Introduction}
Bioconvection is the process through which a spontaneous pattern arises in suspensions of swimming microorganisms. \citep{ref1} The term "bioconvection" was coined by ~\citet{ref3}, while it was ~\citet{ref2} that laid the groundwork for the field with his pioneering studies. Microorganisms have a density that is just slightly higher than the density of the water in which they swim; as a consequence of this, they have a tendency to swim in the opposite direction of the flow of the water, which is upward. In addition, whenever the microbes stop swimming, the patterns disappear along with them. On the other hand, there are also instances of pattern formation that do not need the individuals to up swim or swim at a larger density. \citep{ref1} The types of patterns that are created are determined by a number of different aspects, such as the thickness of the suspension, the number of microorganisms present, and the amount to which they may move about. It has been shown that certain species of flagellated green algae, including \textit{Euglena}, \textit{Chlamydomonas}, \textit{Volvox}, and \textit{Dunaliella}, may generate patterns while they are floating in the water. \citep{ref2,ref4,ref5,ref6} Taxes are the collective name for the reactions that microorganisms display in response to stimuli. These reactions consist of an average swimming motion in a certain direction, and they are exhibited in response to stimuli. There are a wide variety of taxis; however, some of the more known ones include gravitaxis, gyrotaxis, phototaxis, and chemotaxis. There are also many other types of taxis. The response that the microorganisms have to gravity is referred to as gravitaxis. Chemotaxis is a kind of swimming behavior that takes place as a response to chemical gradients, whereas gyrotaxis is created by achieving a balance between the torques that are produced as a consequence of gravity and shear. The movement away from (towards) the direction in which the light source is directed is an expression of negative phototaxis (positive phototaxis). Just phototaxis will be covered in this particular work. Studies have demonstrated that the pattern spacing and size of patterns that are formed by bioconvection may be affected by the intensity of the various kinds of illumination as well as the magnitude of the intensity. \citep{ref2,ref5} The pattern of bioconvection can also alter in terms of its shape, size, and scale depending on the amount of light that is present in the environment. \citep{ref7} The bioconvection pattern shifts for two main reasons: first, the microorganisms will swim toward (or away from) the light source if the total intensity $\mathcal{G}$ is less (or greater) than the critical total intensity $\mathcal{G}_c$.  Microorganisms absorb and scatter light, which is the second reason. \cite{ref8}

The phototaxis model that was suggested by ~\citet{ref15} will be utilized in this study. The Navier-Stokes equation, the radiative transfer equation for controlling the propagation of oblique irradiation, and the cell conservation equation are all used to simulate an incompressible fluid. The light source used to illuminate the algal solution differs significantly from their conceptualization. Due to the fact that the sun hits the surface at a variety of off-normal angles, the current study makes the assumption that the irradiation is collimated in an oblique direction. As a result, the light intensity profiles may be reshuffled by the radiation field throughout the algal suspension, therefore regulating photosynthesis via phototaxis. Since many motile algae must produce their own food through photosynthesis, they exhibit strong phototaxis. Due to the importance of photosynthetic swimming for many motile algae, it is important to consider the effects of forward scattering in realistic and credible models on phototaxis. \cite{gil} The fundamental steady state for a suspension illuminated from above with a finite depth is one in which phototaxis due to diffusion, absorption, and scattering are in equilibrium. This results in the formation of a horizontal sublayer composed of dense microorganisms. The region below a sublayer is gravitationally unstable whereas the region above is gravitationally stable. Critical total intensity $\mathcal{G}_c$ influences the position of the sublayer. The sublayer position is at the top (or bottom) of the suspension if the total intensity $\mathcal{G}$ at every position in the whole suspension is less (or greater) than the critical total intensity $\mathcal{G}_c$. If the critical total intensity $\mathcal{G}_c$ is the same as the total intensity $\mathcal{G}$ inside the suspension, the sublayer position also exists between the top and bottom boundaries. Thus, if one fluid layer becomes unstable, its fluid movements will enter the stable layer in penetrative convection. \citep{ref10}  

~\citet{ref9} initially introduced the concept of phototactic bioconvection. The authors investigated a model that includes vertically collimated irradiation on the non-scattering absorbing suspension. This model of phototaxis and shading is utilized to analyze the linear stability of a suspension of phototactic microorganisms that is evenly irradiated from above and swims in a fluid that has a density that is somewhat lower than that of the algae themselves. After that, ~\citet{ref13} developed a new phototactic model in two dimensions and analyzed its linear stability. A conservative finite-difference approach was employed to numerically solve the governing equations in this model. The impact of scattering was disregarded in these models. The phototactic bioconvection model developed by ~\citet{ref14} incorporated the impact of isotropic scattering. In the equilibrium state for particular parameter values, the microorganisms gather in two horizontal layers at various depths as a consequence of scattering. \citet{kumar2023} extended the work of \citet{ref14} and demonstrated that the system is more stable for the rigid upper surface than the stress-free upper surface. ~\citet{ref15} studied the effect of forward anisotropic scattering on phototactic bioconvection suspension. A detailed computational investigation of the linear stability was provided by the authors, with the focus being placed specifically on the forward scattering impact. ~\citet{ref16} illustrated the influence of both diffused and collimated irradiation on a scattering suspension. A bioconvection model of phototactic irradiation was simulated in \citet{ref16a}, using both diffuse and collimated irradiation on an anisotropic forward scattering suspension.
According to the findings of the author, the base concentration profile of the unimodal moved towards the bimodal and vice versa. ~\citet{ref17} have recently studied a model of phototactic bioconvection in non-scattering suspension with oblique illumination. For specific parameters and angle of incidence variations, the authors observed the most unstable mode transitions from stable to overstable or vice versa. More recently, ~\citet{ref18} has studied oblique irradiation's impact on isotropic scattering suspension. Variation in the angle of incidence revealed two kinds of nature for lower and higher scattering albedo in the fundamental state concentration profile. Nevertheless, there has been no research done on the beginning stages of phototactic bioconvection that takes into account the effects of oblique collimated irradiation on a suspension of algae that exhibits forward anisotropic scattering.
Therefore, research is being conducted to determine how oblique collimated irradiation affects forward anisotropic scattering bioconvection.


\begin{figure*}
    \centering
    \includegraphics[width=18cm, height=13cm]{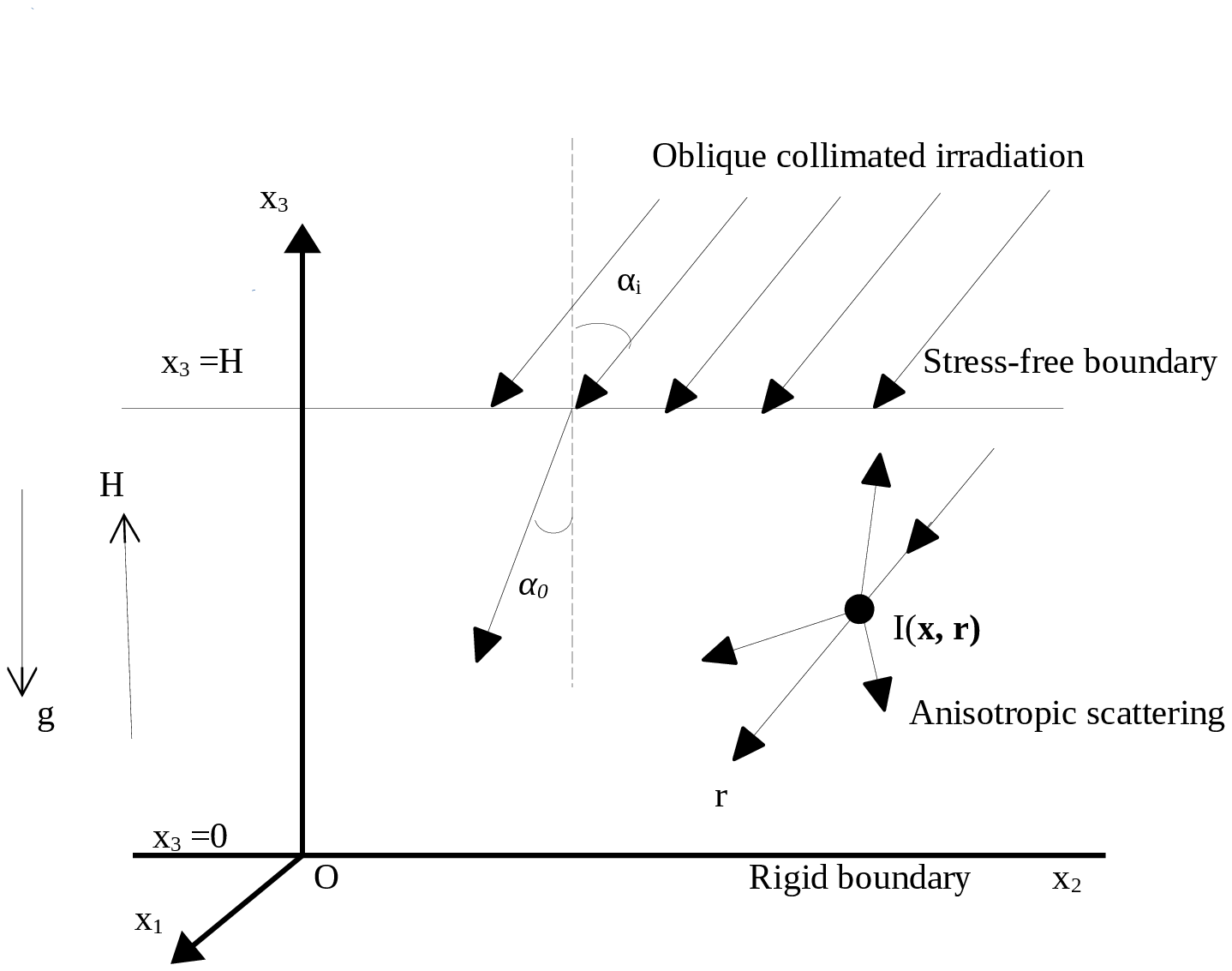}
    \caption{Problem’s geometrical configuration.}
   \label{anisofig.eps}
 \end{figure*}

\section{MATHEMATICAL FORMULATION}

Take into account the motion of phototactic microorganisms suspended within a layer of depth $H$. The lower boundary at $x_3 = 0$ is assumed to be rigid, whereas the upper boundary at $x_3 = H$ is assumed to be stress-free. Let $I(\boldsymbol{x},\boldsymbol{r})$ denote the intensity of the radiation at location $\boldsymbol{x}$ along the unit direction $\boldsymbol{r}=\cos{\alpha}\boldsymbol{k}+\sin{\alpha}(\cos{\zeta}\boldsymbol{i}+\sin{\zeta}\boldsymbol{j})$, where $\boldsymbol{x}$ is a coordinate measured in relation to the rectangular Cartesian axes $Ox_1 x_2 x_3$, with the $x_3$-axis pointing upward (see Figure 1). Also, $\alpha$ and $\zeta$ represent the polar and azimuthal angles of the direction $\boldsymbol{r}$. The horizontal top layer of the phototactic suspension is illuminated by the oblique collimated irradiation from the light source of uniform intensity $I^0$, and the impact of diffuse irradiation is ignored. We are able to calculate the angle of refraction $\alpha_0$ by applying Snell's equation of refraction, which states that ($\sin{\alpha_i} = N_0 \sin{\alpha_0}$). In this expression, $\alpha_i$ indicates the angle of the incidence, and $N_0 = 1.333$ represents the refractive index of water.  Similar to other bioconvection models, we consider a continuous distribution of cells rather than a discrete population of cells.\cite{ref1}  The algal solution is diluted to such an extent that the volume fraction of the cells is reduced to a negligible level, thereby neglecting cell-cell interactions. The cell concentration is denoted by $n$, and $\boldsymbol{v}$ represents the average velocity of all of the material contained within a small volume $\delta V$. The equation of continuity for an incompressible fluid is given by
\begin{equation}
\label{eqn:1}
    \boldsymbol{\nabla} \cdot\boldsymbol{v}=0.
\end{equation}
As negatively buoyant cells disproportionately impact the suspension, we assume that the bulk stress is dominated by Stokeslets and that all other contributions from the cells may be safely disregarded. 
Hence, if we disregard all of the forces acting on the fluid, with the exception of the cell's negative buoyancy, which is represented by the equation $n\vartheta \Delta \varrho g$ per unit volume, where $g$ is the acceleration caused by gravity, the momentum equation according to the Boussinesq approximation is as follows: \cite{ref12}
\begin{equation}
\label{eqn:2}
    \varrho\left(\frac{\partial}{\partial t}+\boldsymbol{v}\cdot\boldsymbol{\nabla}\right)\boldsymbol{v}=\mu\nabla^2\boldsymbol{v}-\boldsymbol{\nabla} \mathcal{P}-n\vartheta \Delta \varrho g \boldsymbol{k},
\end{equation}
 where $t$ is time, $\varrho$ is the density of the fluid, each algal cell has a density $\rho+\Delta \rho$ ($\Delta \rho \ll \rho)$ and the volume $\vartheta$, $\mathcal{P}$ is the excess pressure above hydrostatic, $\mu$ is the dynamic viscosity of the suspension, $\frac{\partial}{\partial t}+\boldsymbol{v}\cdot\boldsymbol{\nabla}=\frac{D}{Dt}$ is the material time derivative, and $\boldsymbol{k}$ is the unit vector in $x_3$ direction. 

The cell conservation equation is given by \cite{ref1,ref9}
\begin{equation}
   \label{eqn:3}
    \frac{\partial n}{\partial t}=-\boldsymbol{\nabla}\cdot\boldsymbol{A}.
\end{equation}
The flux of cells $\boldsymbol{A}$ is given by
\begin{equation}
    \label{eqn:4}
    \boldsymbol{A}=n \boldsymbol{U}_c+n\boldsymbol{v}-\boldsymbol{D}\boldsymbol{\nabla} n.
\end{equation}
The average swimming of the cells $\boldsymbol{U}_c$ generates the first component on the right-hand side of Eq. (\ref{eqn:4}). The second term reflects the flux caused by the cells' advection by the bulk fluid flow, while the third term stands in for the random component of the cells' motion. The diffusivity tensor $\boldsymbol{D}$ is chosen to be isotropic and constant, thus $\boldsymbol{D}= D I$. Here, $D$ and $I$ are the diffusion coefficient and identity tensor, respectively. We assume that individual algal cells are identical spheres whose centers of gravity and buoyancy are both located at the same point in space and that they are exclusively phototactic. Throughout the process of calculating the cell flux vector for the proposed model, two main presumptions have been utilized. First, the influence of viscous torque owing to shear in the flow, which can contribute to the horizontal component of the mean swimming direction, is disregarded since the algal cells are exclusively phototactic. Second, the diffusion tensor, which ought to be a function of the intensity of the light, is considered to be a constant isotropic tensor rather than being derived from a swimming velocity auto-correlation function. This is done in order to simplify the analysis. So, this model may be assumed to be true in the limited scenario in order to comprehend the complexity of bioconvection owing to phototaxis before examining more comprehensive complicated models. Based on these hypotheses, we are able to remove the Fokker-Planck equation from the governing equations of phototactic bioconvection. 

For an absorbing and scattering medium, the radiative transfer equation, which will be abbreviated as RTE from here on out, is defined as \cite{ref-chand,ref-modest}
\begin{equation}
\label{eqn:5}
    \boldsymbol{r}\cdot\boldsymbol{\nabla}I(\boldsymbol{x},\boldsymbol{r})+(\kappa+\sigma)I(\boldsymbol{x},\boldsymbol{r})=\frac{\sigma}{4\pi}\int_0^{4\pi}I(\boldsymbol{x},\boldsymbol{r}')\Psi(\boldsymbol{r}',\boldsymbol{r})d\Omega',
\end{equation}
where $\kappa$, $\sigma$, and $\Omega'$ are the absorption coefficient, scattering coefficient, and solid angles, respectively. Scattering from the direction $\boldsymbol{r}'$ to the direction $\boldsymbol{r}$ is represented by the probability density function, $\Psi(\boldsymbol{r}',\boldsymbol{r})$, referred to as the scattering phase function. For the reason of convenience, we consider that $\Psi(\boldsymbol{r}',\boldsymbol{r})$ is linearly anisotropic and has azimuthal symmetry: \cite{ref-modest}
\begin{equation}
\label{eqn:6}
    \Psi(\boldsymbol{r}',\boldsymbol{r})=1+A \cos{\alpha}\cos{\alpha'}.
\end{equation}
The unit vectors $\boldsymbol{r}$ and $\boldsymbol{r}'$ are described by the polar angles $\alpha$ and $\alpha'$ from the $x_3$-axis, respectively. 
The linearly anisotropic scattering coefficient is denoted by $A\in [-1,1]$. The scattering happens in a forward direction for $A>0$ and in a backward direction for $A<0$, respectively. Isotropic scattering is represented by the scenario in which $A = 0$. The intensity of the radiation that is emitted from the top boundary surface at the positions $\boldsymbol{x}_b=(x_1,x_2, H)$ is given by: 
\begin{equation}
\label{eqn:7}
   I(\boldsymbol{x}_b, \boldsymbol{r})=\delta(\boldsymbol{r}-\boldsymbol{r}^0)=I^0\delta(\cos{\alpha}-\cos{\alpha_0})\delta(\zeta-\zeta_0).
\end{equation}
In this equation, $I^0$ represents the incidence collimated radiation, $\boldsymbol{r}^0=\cos{\alpha_0}\boldsymbol{k}+\sin{\alpha_0}(\cos{\zeta_0}\boldsymbol{i}+\sin{\zeta_0}\boldsymbol{j})$ represents the incident direction, and $\zeta$ represents the azimuthal angle of the direction. Unit vectors along the axes $x_1$, $x_2$, and $x_3$ are denoted by the symbols $\boldsymbol{i}$, $\boldsymbol{j}$, and $\boldsymbol{k}$, respectively. 
The Dirac-delta function $\delta$ satisfies \cite{ref-modest}
\begin{equation*}
    \int_0^{4\pi} f(\boldsymbol{r})\delta(\boldsymbol{r}-\boldsymbol{r}^0) d\Omega=f(\boldsymbol{r}^0)=f(\eta_3^0,\zeta_0), \quad \eta_3^0=\cos{\alpha_0}.
\end{equation*}
The concentration $n$ has a linear relation with the absorption and scattering coefficients so that $\kappa=\varpi n$ and $\sigma=\varsigma n$. Regarding the scattering albedo $\omega=\sigma/(\kappa+\sigma)\in[0,1]$, the RTE in a linearly anisotropic scattering algal suspension hence becomes
\begin{equation}
\label{eqn:8}
    \boldsymbol{r}\cdot\boldsymbol{\nabla}I(\boldsymbol{x},\boldsymbol{r})+\beta I(\boldsymbol{x},\boldsymbol{r})
    =\frac{\omega \beta}{4\pi}\int_0^{4\pi}I(\boldsymbol{x},\boldsymbol{r}')(1+A \cos{\alpha}\cos{\alpha'})d\Omega',
\end{equation}
where $\beta=(\varpi+\varsigma)n$ is the extinction coefficient. The total intensity $\mathcal{G}(\boldsymbol{x})$ and radiative heat flux $\boldsymbol{q}(\boldsymbol{x})$ at a point $\boldsymbol{x}$ are given by\cite{ref-modest}
\begin{equation}
\label{eqn:9}
    \mathcal{G}(\boldsymbol{x})=\int_0^{4\pi}I(\boldsymbol{x},\boldsymbol{r})d\Omega,
\end{equation}
\begin{equation}
\label{eqn:10}
    \boldsymbol{q}(\boldsymbol{x})=\int_0^{4\pi}I(\boldsymbol{x},\boldsymbol{r})\boldsymbol{r}d\Omega.
\end{equation}
Let $\boldsymbol{P}$ be the swimming direction's unit vector, and $<\boldsymbol{P}>$ be the ensemble average of the swimming direction for all cells in an elemental volume. Many microorganisms' swimming speed is independent of light, position, time, and direction. \cite{hill1997biased} All cells are considered to move through the fluid at the same rate. $U_c$ denotes the ensemble swimming speed. Therefore, average swimming velocity is given by
\begin{equation}
\label{eqn:11}
    \boldsymbol{U}_c=U_c<\boldsymbol{P}>.
\end{equation}
The mean swimming direction, $<\boldsymbol{P}>$, is given by
\begin{equation}
\label{eqn:12}
    <\boldsymbol{P}>=-T(\mathcal{G})\frac{\boldsymbol{q}}{\iota+|\boldsymbol{q}|},
\end{equation}
where $\iota\ge0$ is a constant. In the case of collimated irradiation, which is considered throughout this study, the light intensity in the medium is not uniformly distributed, therefore we can use $\iota=0$. Thus, in this study, we consider $\iota=0$. The negative sign represents the fact that a microorganism receives light coming from the opposite direction of the radiative heat flux. Phototaxis function $T(\mathcal{G})$ is defined as
\begin{equation}
\label{eqn:13}
 T(\mathcal{G}) = \left\{ \begin{array}{lll}
         < 0 & \mbox{if $\mathcal{G}_c < \mathcal{G}$};\\
        \ge 0 & \mbox{if $\mathcal{G}_c \ge \mathcal{G}$}.\end{array} \right. 
\end{equation}

It is considered that the lower boundary is rigid, whereas the upper boundary is stress-free. The boundary conditions are
\begin{equation}
\label{eqn:14}
    \boldsymbol{v}=\boldsymbol{v}\times\boldsymbol{k}=\boldsymbol{A}\cdot\boldsymbol{k}=0 \quad \text{at} \quad x_3=0,
\end{equation}
\begin{equation}
\label{eqn:15}
    \boldsymbol{v}\cdot\boldsymbol{k}=\frac{\partial^2}{\partial x_3^2}(\boldsymbol{v}\cdot\boldsymbol{k})=\boldsymbol{A}\cdot\boldsymbol{k}=0 \quad \text{at} \quad x_3=H.
\end{equation}
Radiation traveling vertically and perpendicular to the $x_3$-axis is directed at the top boundary. It is assumed that the top and bottom bounds are non-reflective, such that 
\begin{equation}
\label{eqn:16}
    I(x_1,x_2,H,\alpha,\zeta)=I^0\delta(\cos{\alpha}-\cos{\alpha_0})\delta(\zeta-\zeta_0),\quad
    \pi/2\le\alpha\le\pi.
\end{equation}
\begin{equation}
\label{eqn:17}
    I(x_1,x_2,0,\alpha,\zeta)=0,\quad 0\le\alpha\le\pi/2,
\end{equation}

All lengths, cell concentration, time, fluid velocity, and pressure are scaled by $H$, $\bar{n}$, $H^{2}/D$, $D/H$, and $\mu D/H^{2}$ to produce the non-dimensional governing system of bioconvection equations. 

After substituting the dimensionless variables in governing equations
\begin{equation}
\label{eqn:18}
    \boldsymbol{\nabla} \cdot\boldsymbol{v}=0,
\end{equation}
\begin{equation}
\label{eqn:18}
    \frac{1}{S_c}\left(\frac{\partial}{\partial t}+\boldsymbol{v}\cdot\boldsymbol{\nabla}\right)\boldsymbol{v}=\nabla^2\boldsymbol{v}-\boldsymbol{\nabla} \mathcal{P}-R_a n\boldsymbol{k},
\end{equation}
\begin{equation}
\label{eqn:20}
    \frac{\partial n}{\partial t}=-\boldsymbol{\nabla}\cdot(n U_s\boldsymbol{P}+n\boldsymbol{v}-\boldsymbol{\nabla} n).
\end{equation}

Here $S_c=\frac{\nu}{D}$ is the Schmidt number, $U_s=\frac{U_c H}{D}$ is the dimensionless swimming speed, $R_a=\frac{\bar{n} \vartheta \Delta \rho g H^3}{\mu D}$ is the Rayleigh number, $\nu=\mu/\rho$ is the kinematic viscosity.

Dimensionless boundary conditions for rigid surfaces become
\begin{equation}
\label{eqn:equation21}
    \boldsymbol{v}=\boldsymbol{v}\times\boldsymbol{k}=(n U_s\boldsymbol{P}+n\boldsymbol{v}-\boldsymbol{\nabla} n)\cdot\boldsymbol{k}=0 \quad \text{at} \quad x_3=0,
\end{equation}
while for a stress-free surface
\begin{equation}
\label{eqn:equation22}
    \boldsymbol{v}\cdot\boldsymbol{k}=\frac{\partial^2}{\partial x_3^2}(\boldsymbol{v}\cdot\boldsymbol{k})=(n U_s\boldsymbol{P}+n\boldsymbol{v}-\boldsymbol{\nabla} n)\cdot\boldsymbol{k}=0 \quad \text{at} \quad x_3=1.
\end{equation}

The dimensionless RTE is 
\begin{equation}
\label{eqn:23}
    \boldsymbol{r}\cdot\boldsymbol{\nabla}I(\boldsymbol{x},\boldsymbol{r})+\tau_h n I(\boldsymbol{x},\boldsymbol{r})
    =\frac{\omega \tau_h n}{4\pi}\int_0^{4\pi}I(\boldsymbol{x},\boldsymbol{r}')(1+A \cos{\alpha}\cos{\alpha'})d\Omega',
\end{equation}
here  $\tau_h=(\varpi+\varsigma)\bar{n} H$ is the optical depth of the suspension. In terms of the direction cosines, $(\eta_1,\eta_2\eta_3)$ of the unit vector $\boldsymbol{r}$, the non-dimensional RTE can be written as
\begin{equation}
\label{eqn:24}
    \eta_1 \frac{\partial I}{\partial x_1}+\eta_2 \frac{\partial I}{\partial x_2}+\eta_3 \frac{\partial I}{\partial x_3}+\tau_h n I(\boldsymbol{x},\boldsymbol{r})
    =\frac{\omega \tau_h n}{4\pi}\int_0^{4\pi}I(\boldsymbol{x},\boldsymbol{r}')(1+A \cos{\alpha}\cos{\alpha'})d\Omega'.
\end{equation}
For their non-dimensional counterparts, the symbols $\mathcal{G}$ and $\boldsymbol{q}$ remain the same. The exact functional form of $T (\mathcal{G})$ depends on the species of microorganisms. \cite{ref9}
For instance, the following is an example of a typical phototaxis function, which is defined in mathematical form: 
\begin{eqnarray}
\label{eqn:24a}
    T(\mathcal{G})=0.8\sin{(1.5\pi\varphi(\mathcal{G}))}-0.1\sin{(0.5\pi\varphi(\mathcal{G}))}, \nonumber \\
    \varphi(\mathcal{G})=0.4 \mathcal{G}\exp{(\Upsilon(2.5-\mathcal{G}))}.
\end{eqnarray}
Parameter $\Upsilon$ is connected to the determination of the value of the critical total intensity. 

In a non-dimensional representation, the top and bottom intensities become 
\begin{equation}
\label{eqn:25}
    I(x_1,x_2,1,\alpha,\zeta)=I^0\delta(\boldsymbol{r}-\boldsymbol{r}^0) , \quad \pi/2\le\alpha\le\pi,
\end{equation}
\begin{equation}
\label{eqn:26}
    I(x_1,x_2,0,\alpha,\zeta)=0,\quad 0\le\alpha\le\pi/2.
\end{equation}

\section{Basic state solution}

Equations (\ref{eqn:18})–(\ref{eqn:20}) and (\ref{eqn:24}), in conjunction with the boundary conditions, have a solution for the static equilibrium that includes: $\boldsymbol{v}=0$, $n=n_b(x_3)$ and $I(x_3,\alpha)$.

The total intensity and the radiative heat flux at steady-state are given by
\begin{equation}
\label{eqn:27}
    \mathcal{G}_b(x_3)=\int_0^{4\pi} I_b(x_3,\alpha)d\Omega,
\end{equation}
\begin{equation}
\label{eqn:28}
   \boldsymbol{q}_b(x_3)=\int_0^{4\pi} I_b(x_3,\alpha)\boldsymbol{r}d\Omega.
\end{equation}
Since $I_b(x_3,\alpha)$ is independent of $\zeta$, the $x_1$ and $x_2$ components of $\boldsymbol{q}_b$ vanish. Therefore,$\boldsymbol{q}_b=q_b\boldsymbol{k}$,  where $q_b=|\boldsymbol{q}_b|$.
The radiative transfer equation at steady-state becomes
\begin{equation}
\label{eqn:29}
    \frac{\partial I_b}{\partial x_3}+\frac{\tau_h n_b I_b}{\eta_3}=\frac{\omega \tau_h n_b}{4\pi \eta_3}(\mathcal{G}_b(x_3)-A q_b \eta_3).
\end{equation}
We decompose steady-state intensity into collimated $I_b^c$ and diffuse parts $I_b^d$, i.e., $I_b=I_b^c+I_b^d$. The collimated component $I_b^c$ satisfies
\begin{equation}
\label{eqn:30}
    \frac{d I_b^c}{d x_3}+\frac{\tau_h n_b I_b^c}{\eta_3}=0,
\end{equation}
with boundary condition
\begin{equation}
\label{eqn:31}
    I_b^c(x_3,\alpha) =I^0\delta(\boldsymbol{r}-\boldsymbol{r}^0), \quad \textrm{at}, \quad z=1.
\end{equation}
Solving Eqs. (\ref{eqn:30}) and (\ref{eqn:31}), we get
\begin{equation}
\label{eqn:32}
    I_b^c=I^0\exp{\left(-\int_1^{x_3}\frac{\tau_h n_b(x_3')}{\eta_3}d x_3'\right)}\delta(\boldsymbol{r}-\boldsymbol{r}^0).
\end{equation}
Now the diffused component $I_b^d$ satisfies
\begin{equation}
\label{eqn:33}
    \frac{d I_b^d}{d x_3}+\frac{\tau_h n_b I_b^d}{\eta_3}=\frac{\omega \tau_h n_b}{4\pi \eta_3}(\mathcal{G}_b(x_3)-A q_b \eta_3),
\end{equation}
with boundary conditions
\begin{equation}
\label{eqn:34}
    I_b^d(1,\alpha)=0, \quad \pi/2\le\alpha\le\pi,
\end{equation}
\begin{equation}
\label{eqn:35}
    I_b^d(0,\alpha)=0,\quad 0\le\alpha\le\pi/2.
\end{equation}
We define a new variable $\tau=\int_{x_3}^1 \tau_h n_b(x_3')d x_3'$. Now, the total intensity $\mathcal{G}_b$ and radiative heat flux $\boldsymbol{q}_b$ become functions of $\tau$ only. In terms of the variable $\tau$, Eq. (\ref{eqn:33}) for the outgoing $(0<\alpha<\pi/2)$ radiation intensity, $I_b^{d+}(\tau,\eta_3)$, is written as
\begin{equation}
\label{eqn:36}
    \eta_3\frac{d I_b^{d+}}{d \tau}-I_b^{d+}(\tau,\eta_3)=-\frac{\omega }{4\pi}(\mathcal{G}_b(\tau)-A q_b(\tau) \eta_3),\quad
    0<\eta_3<1,
\end{equation}
with boundary condition
\begin{equation}
\label{eqn:37}
    I_b^{d+}(\tau_h,\eta_3)=0, \quad 0<\eta_3<1. 
\end{equation}
On solving Eqs. (\ref{eqn:36})-(\ref{eqn:37})),
\begin{equation}
\label{eqn:38}
    I_b^{d+}(\tau,\eta_3)=-\frac{\omega}{4 \pi \eta_3}\int_{\tau_h}^{\tau}(\mathcal{G}_b(\tau')-A q_b(\tau') \eta_3)\exp{\left(\frac{\tau-\tau'}{\eta_3}\right)}d \tau', \quad
    0<\eta_3<1.
\end{equation}
Similarly, the incoming $(\pi/2<\alpha<\pi)$ radiation intensity, $I_b^{d-}(\tau,\eta_3)$, is written as
\begin{equation}
\label{eqn:39}
    I_b^{d-}(\tau,\eta_3)=-\frac{\omega}{4 \pi \eta_3}\int_{0}^{\tau}(\mathcal{G}_b(\tau')-A q_b(\tau') \eta_3)\exp{\left(\frac{\tau-\tau'}{\eta_3}\right)}d \tau', \quad
    -1<\eta_3<0.
\end{equation}
The basic total intensity and radiative heat flux are written as
\begin{equation}
\label{eqn:40}
    \mathcal{G}_b=\mathcal{G}_b^c+\mathcal{G}_b^d,
\end{equation}
\begin{equation}
\label{eqn:41}
    \boldsymbol{q}_b=\boldsymbol{q}_b^c+\boldsymbol{q}_b^d,
\end{equation}
where
\begin{equation*}
    \mathcal{G}_b^c(x_3)=\int_0^{4\pi} I_b^c(x_3,\alpha)d\Omega=I^0\exp{\left(\frac{\tau_h}{\cos{\alpha_0}}\int_1^{x_3}n_b(x_3')d x_3'\right)},
\end{equation*}
\begin{equation*}
    \mathcal{G}_b^d(x_3)=\int_0^{4\pi} I_b^d(x_3,\alpha)d\Omega,
\end{equation*}
\begin{equation*}
    \boldsymbol{q}_b^c(x_3)=\int_0^{4\pi} I_b^c(x_3,\alpha)\boldsymbol{r}d\Omega=-I^0\cos{\alpha_0}\exp{\left(\frac{\tau_h}{\cos{\alpha_0}}\int_1^{x_3}n_b(x_3')d x_3'\right)}\boldsymbol{k},
\end{equation*}
\begin{equation*}
    \boldsymbol{q}_b^d(x_3)=\int_0^{4\pi} I_b^d(x_3,\alpha)\boldsymbol{r}d\Omega.
\end{equation*}
Substitution of Eqs. (\ref{eqn:38}) and (\ref{eqn:39}) into Eqs. (\ref{eqn:40}) and (\ref{eqn:41}), yield two coupled Fredholm integral equations of the second kind\cite{ref-modest,sarma2005}:
\begin{equation}
\label{eqn:42}
    \mathcal{G}_b(\tau)=e^{-\frac{\tau}{\cos{\alpha_0}}}+\frac{\omega}{2}\int_0^{\tau_h}(\mathcal{G}_b(\tau')E1(|\tau-\tau'|)+
    sgn(\tau-\tau')A q_b(\tau')E2(|\tau-\tau'|)d \tau',
\end{equation}
\begin{equation}
\label{eqn:43}
    q_b(\tau)=e^{-\frac{\tau}{\cos{\alpha_0}}}+\frac{\omega}{2}\int_0^{\tau_h}(A q_b(\tau')E3(|\tau-\tau'|)+
    sgn(\tau-\tau') \mathcal{G}_b(\tau')E2(|\tau-\tau'|)d \tau'.
\end{equation}
In this context, 'sgn' refers to the sign function, whereas $E1 (x_1)$, $E2 (x_2)$, and $E3 (x_3)$ are the first, second, and third-order exponential integral functions, respectively. \cite{ref-chand} Equations (\ref{eqn:42}) and (\ref{eqn:43}) can be solved by employing a technique known as the subtraction of singularity. \cite{numerical-res} In the case of isotropic scattering, it is important to note that these two equations become uncoupled from one another. \cite{ref18}

The mean swimming direction in the basic state becomes
\begin{equation*}
    <\boldsymbol{P}_b>=-T_b\frac{\boldsymbol{q}_b}{|\boldsymbol{q}_b|}=T_b\boldsymbol{k},
\end{equation*}
where
\begin{equation*}
    T_b=T(\mathcal{G}_b).
\end{equation*}
The steady-state cell conservation equation is written as 
\begin{equation}
\label{eqn:44}
\frac{d n_b}{dx_3}-U_s T_b n_b=0,
\end{equation}
which is supplemented by the cell conservation relation
\begin{equation}
\label{eqn:45}
    \int_0^1 n_b(x_3)dx_3=1.
\end{equation}
The boundary value problem formed by equations (28)-(31) is resolved numerically by the shooting technique.

\section{Linear and normal modes analysis}
A small perturbation $\epsilon(0<\epsilon \ll 1)$ is made in the basic state to examine the linear instability. 
$\boldsymbol{v}=\boldsymbol{0}+\epsilon \boldsymbol{v}^*(x_1,x_2,x_3,t)+O(\epsilon^2)$, $n=n_b(x_3)+\epsilon n^*(x_1,x_2,x_3,t)+O(\epsilon^2)$, $\mathcal{P}=\mathcal{P}_b+\epsilon \mathcal{P}^*+O(\epsilon^2)$, $<\boldsymbol{P}>=<\boldsymbol{P}_b>+\epsilon<\boldsymbol{P}^*>+O(\epsilon^2)$, $\mathcal{G}=\mathcal{G}_b+\epsilon\mathcal{G}^*+O(\epsilon^2)$, $\boldsymbol{q}=\boldsymbol{q}_b+\boldsymbol{q}^*+O(\epsilon^2)$, $I=(I_b^c+I_b^d)+\epsilon (I^{*c}+I^{*d})+O(\epsilon^2)$, where $\boldsymbol{v}^* = (u^* , v^* , w^*)$.
Thus, the linearized governing equations are written as
\begin{equation}
\label{eqn:46}
    \nabla \cdot \boldsymbol{v}^* =0  
\end{equation}
\begin{equation}
\label{eqn:47}
   \frac{1}{S_c}\frac{\partial \boldsymbol{v}^*}{\partial t}=\nabla^2\boldsymbol{v}^*-\boldsymbol{\nabla} \mathcal{P}^*-n^* R_a \boldsymbol{k},
\end{equation}

\begin{equation}
\label{eqn:48}
    \frac{\partial n^*}{\partial t}+\frac{dn_b}{dx_3}w^*=\nabla^2 n^*-U_s\boldsymbol{\nabla}\cdot(n^*<\boldsymbol{P}_b>+n_b<\boldsymbol{P}^*>).
\end{equation}

The total intensity $\mathcal{G}$ and radiative heat flux $\boldsymbol{q}$ can be written as $\mathcal{G}=\mathcal{G}_b+\epsilon\mathcal{G^*}+O(\epsilon^2)=(\mathcal{G}_b^c+\mathcal{G}_b^d)+\epsilon(\mathcal{G}^{*c}+\mathcal{G}^{* d})+O(\epsilon^2)$, $\boldsymbol{q}=\boldsymbol{q}_b+\epsilon\boldsymbol{q}^*+O(\epsilon^2)=(\boldsymbol{q}_b^c+\boldsymbol{q}_b^d)+\epsilon(\boldsymbol{q}^{*c}+\boldsymbol{q}^{* d})+O(\epsilon^2)$. Thus, the steady collimated total intensity is perturbed as 
\begin{equation*}
    \mathcal{G}_b^c+\epsilon\mathcal{G}^{*c}+O(\epsilon^2)=
    I^0\exp{\left(\frac{\tau_h}{\cos{\alpha_0}}\int_1^{x_3}(n_b(x_3')+\epsilon n^*+O(\epsilon^2)d x_3'\right)},
\end{equation*}
On collecting $O(\epsilon)$ terms
\begin{equation}
\label{eqn:49}
    \mathcal{G}^{*c}=I^0\exp{\left(\frac{\tau_h}{\cos{\alpha}_0}\int_1^{x_3}n_b(x_3')d x_3'\right)}\left(\frac{\tau_h}{\cos{\alpha}_0}\int_1^{x_3}n^* d x_3'\right).
\end{equation}
Similarly, 
\begin{equation}
\label{eqn:50}
    \mathcal{G}^{*d}=\int_0^{4\pi}I^{*d}(\boldsymbol{x},\boldsymbol{r})d \Omega,
\end{equation}
\begin{equation}
\label{eqn:51}
    \boldsymbol{q}^{*c}=-I^0\cos{\alpha}_0\exp{\left(\frac{\tau_h}{\cos{\alpha}_0}\int_1^{x_3}n_b(x_3')d x_3'\right)}
    \left(\frac{\tau_h}{\cos{\alpha}_0}\int_1^{x_3}n^* d x_3'\right)\boldsymbol{k},
\end{equation}
\begin{equation}
\label{eqn:52}
    \boldsymbol{q}^{*d}=\int_0^{4\pi}I^{*d}(\boldsymbol{x},\boldsymbol{r})\boldsymbol{r}d \Omega.
\end{equation}

The average swimming orientation is written as
\begin{equation}
\label{eqn:53}
    <\boldsymbol{P}>=<\boldsymbol{P}_b>+\epsilon<\boldsymbol{P}^*> +O(\epsilon^2)
    =T(\mathcal{G}_b+\epsilon \mathcal{G}^*+O(\epsilon^2))\frac{\boldsymbol{q}_b+\boldsymbol{q}^*+O(\epsilon^2)}{|\boldsymbol{q}_b+\boldsymbol{q}^*+O(\epsilon^2)|}.
\end{equation}
On collecting $O(\epsilon)$ terms
\begin{equation}
\label{eqn:54}
    <\boldsymbol{P}^*>=\mathcal{G}^*\frac{\partial T_b}{\partial \mathcal{G}}\boldsymbol{k}-T_b\frac{\boldsymbol{q_H^*}}{q_b},
\end{equation}
where $\boldsymbol{q_H^*}$ is the horizontal component of the perturbed intensity flux $\boldsymbol{q}^*$.
Substituting Eq. (\ref{eqn:54}) into Eq. (\ref{eqn:48}) and simplifying we get
\begin{equation}
\label{eqn:55}
    \frac{\partial n^*}{\partial t}+\frac{dn_b}{dx_3}w^*=\nabla^2 n^*-U_s\frac{\partial}{\partial x_3}\left(T_b n^*+n_b\frac{d T_b}{d \mathcal{G}}\mathcal{G}^*\right)+
    U_s n_b\frac{T_b}{q_b}\left(\frac{\partial q_1^*}{\partial x_1}+\frac{\partial q_2^*}{\partial x_2}\right).
\end{equation}
Taking the curl of Eq. (\ref{eqn:47}) twice and getting the $x_3$-component eliminates $\mathcal{P}^*$ and the horizontal component of $\boldsymbol{v}$. Hence, Equations (\ref{eqn:46}), (\ref{eqn:47}), and (\ref{eqn:55}) may be simplified into two equations for $w^*$ and $n^*$. These quantities can then be decomposed into their normal modes in such a way that 
\begin{equation*}
    w^*=\hat{w}(x_3)\exp{[\gamma t+i(a_1 x+a_2 y)]},
\end{equation*}
\begin{equation*}
    n^*=\hat{n}(x_3)\exp{[\gamma t+i(a_1 x+a_2 y)]},
\end{equation*}
Here, $a_1$ and $a_2$ are wavenumbers in the $x_1$ and $x_2$ directions and the resultant $a=\sqrt{a_1^2+a_2^2}$ is a horizontal wavenumber. The growth rate is represented by $Re(\gamma)$.
From Eq. (\ref{eqn:24}), the perturbed diffused radiation intensity $I^{*d}$ satisfies
\begin{equation}
\label{eqn:56}
    \eta_1 \frac{\partial I^{*d}}{\partial x_1}+\eta_2 \frac{\partial I^{*d}}{\partial x_2}+\eta_3 \frac{\partial I^{*d}}{\partial x_3}+\tau_h n_b I^{*d}(\boldsymbol{x},\boldsymbol{r})
    =\frac{\omega \tau_h }{4\pi}(n_b \mathcal{G}^* +\mathcal{G}_b n^*+A \eta_3(n_b \boldsymbol{q}^*\cdot\boldsymbol{k}-q_b n^*))-\tau_h I_b n^*,
\end{equation}
with boundary conditions
\begin{equation}
\label{eqn:57}
    I^{*d}(x_1,x_2,1,\eta_1,\eta_2,\eta_3)=0,\quad \pi/2\le \alpha\le\pi,\quad 0\le\zeta\le 2\pi,
\end{equation}
\begin{equation}
\label{eqn:58}
    I^{*d}(x_1,x_2,0,\eta_1,\eta_2,\eta_3)=0,\quad 0\le \alpha\le\pi/2,\quad 0\le\zeta\le 2\pi.
\end{equation}
In normal modes, $I^{*d}$ can be written as
\begin{equation*}
    I^{*d}=\psi^d(x_3,\eta_1,\eta_2,\eta_3)\exp{[\gamma t+i(a_1 x+a_2 y)]},
\end{equation*}
From Equations (\ref{eqn:49}) and (\ref{eqn:50}) we get
\begin{equation}
\label{eqn:59}
    \mathcal{G}^{*c}=G^c(x_3)\exp{[\gamma t+i(a_1 x+a_2 y)]},
\end{equation}
\begin{equation}
\label{eqn:60}
    \mathcal{G}^{*d}=G^d(x_3)\exp{[\gamma t+i(a_1 x+a_2 y)]},
\end{equation}

where

\begin{equation}
\label{eqn:61}
    G^c(x_3)=I^0\exp{\left(\frac{\kappa}{\cos{\alpha}_0}\int_1^{x_3}n_b(x_3')d x_3'\right)}\left(\frac{\kappa}{\cos{\alpha}_0}\int_1^{x_3}\hat{n} d x_3'\right),
\end{equation}
\begin{equation}
\label{eqn:62}
    G^d(x_3)=\int_0^{4\pi}\psi^{*d}(x_3,\eta_1,\eta_2,\eta_3)d \Omega.
\end{equation}
Similarly,
\begin{equation*}
    (q_1^*,q_2^*,q_3^*)=(\hat{q}_1(x_3),\hat{q}_2(x_3),\hat{q}_3(x_3))\exp{[\gamma t+i(a_1 x+a_2 y)]},
\end{equation*}
where
\begin{eqnarray*}
    \hat{q}_1(x_3)=\int_0^{4\pi}\psi^{*d}(x_3,\eta_1,\eta_2,\eta_3)\eta_1 d \Omega, \nonumber \\
    \hat{q}_2(x_3)=\int_0^{4\pi}\psi^{*d}(x_3,\eta_1,\eta_2,\eta_3)\eta_2 d \Omega,\nonumber \\
    \hat{q}_3(x_3)=\int_0^{4\pi}\psi^{*d}(x_3,\eta_1,\eta_2,\eta_3)\eta_3 d \Omega.
\end{eqnarray*}
Now, Equation (\ref{eqn:56}) with boundary conditions (\ref{eqn:57})-(\ref{eqn:58}) can be written as
\begin{equation}
\label{eqn:63}
    \frac{\partial \psi^d}{\partial x_3}+\frac{i(a_1 \eta_1+a_2\eta_2)+\tau_h n_b}{\eta_3}\psi^d=
    \frac{\omega \tau_h}{4\pi \eta_3}(n_b G+\mathcal{G}_b \hat{n}+A \eta_3(n_b\hat{q}_3-q_b\hat{n}))-\frac{\tau_h}{\eta_3}I_b^d\hat{n},
\end{equation}
boundary conditions
\begin{equation}
\label{eqn:64}
    \psi^{d}(1,\eta_1,\eta_2,\eta_3)=0,\quad \pi/2\le \alpha\le\pi,\quad 0\le\zeta\le 2\pi,
\end{equation}
\begin{equation}
\label{eqn:65}
    \psi^{d}(0,\eta_1,\eta_2,\eta_3)=0,\quad 0\le \alpha\le\pi/2,\quad 0\le\zeta\le 2\pi.
\end{equation}

The governing equations in normal modes become 
\begin{equation}
    \label{eqn:66}
    \gamma S_c^{-1}\left(\frac{d^2}{dx_3^2}-a^2\right) \hat{w}(x_3) -\left(\frac{d^2}{dx_3^2}-a^2\right)^2\hat{w}(x_3)
    = a^2R_a\hat{n}(x_3),
\end{equation}

\begin{equation}
\label{eqn:67}
    \left(\gamma+a^2-\frac{d^2}{d x_3^2}\right)\hat{n}+U_s\frac{d}{d x_3}\left(T_b\hat{n}+n_b\frac{d T_b}{d \mathcal{G}}G\right)-
    i\frac{U_s n_b T_b}{q_b}(a_1\hat{q}_1+a_2\hat{q}_2)=-\frac{d n_b}{d x_3}\hat{w},
\end{equation}

Boundary conditions in normal modes
\begin{equation}
    \label{eqn:68}
    \hat{w}=\frac{d\hat{w}}{d x_3}=\frac{d \hat{n}}{d x_3}-U_s T_b \hat{n}-U_s n_b\frac{d T_b}{d \mathcal{G}}G=0, \quad x_3=0,
\end{equation}
\begin{equation}
    \label{eqn:69}
    \hat{w}=\frac{d^2\hat{w}}{d x_3^2}=\frac{d \hat{n}}{d x_3}-U_s T_b \hat{n}-U_s n_b\frac{d T_b}{d \mathcal{G}}G=0, \quad x_3=1.
\end{equation}
Define a new variable
\begin{equation}
\label{eqn:70}
    N(x_3)=\int_1^{x_3} \hat{n} d \bar{x}_3,
\end{equation}
so that the system of equations becomes
\begin{equation}
\label{eqn:71}
    \frac{d^4 \hat{w}}{d x_3}-(2 a^2+\gamma S_c^{-1})\frac{d^2 \hat{w}}{d x_3^2}+a^2(a^2+\gamma S_c^{-1})\hat{w} 
    =-a^2 R_a \frac{d N}{d x_3},
\end{equation}

\begin{eqnarray}
\label{eqn:72}
    U_s\frac{d}{d x_3}\left(n_b\frac{d T_b}{d \mathcal{G}}G^d\right)-i\frac{U_s n_b T_b}{q_b}(a_1 \hat{q}_1+a_2\hat{q}_2)+
    \frac{\tau_h}{\cos{\alpha}_0}U_s\frac{d}{d x_3}\left(n_b\mathcal{G}_b^c\frac{d T_b}{d \mathcal{G}}\right)N+\nonumber \\
    \left(\gamma+a^2+\frac{2\tau_h}{\cos{\alpha}_0}U_s n_b\mathcal{G}_b^c\frac{d T_b}{d \mathcal{G}}+U_s\frac{d T_b}{d \mathcal{G}}\frac{d\mathcal{G}_b^d}{d x_3}\right)\frac{d N}{d x_3}+
    U_s T_b \frac{d^2 N}{d x_3^2}-\frac{d^3 N}{d x_3^3}=-\frac{d n_b}{d x_3}\hat{w}.
\end{eqnarray}
Also, boundary conditions become
\begin{equation}
    \label{eqn:73}
    \hat{w}=\frac{d\hat{w}}{d x_3}=-U_s T_b \frac{d N}{d x_3}-U_s n_b\frac{d T_b}{d \mathcal{G}}G+\frac{d^2 N}{d x_3^2}=0,
    \end{equation}
\begin{equation}
    \label{eqn:74}
    \hat{w}=\frac{d^2\hat{w}}{d x_3^2}=-U_s T_b \frac{d N}{d x_3}-U_s n_b\frac{d T_b}{d \mathcal{G}}G+\frac{d^2 N}{d x_3^2}=0,
    \end{equation}
and 
\begin{equation}
\label{eqn:75}
    N(x_3)=0 \quad \text{at} \quad x_3=1.
\end{equation}

\begin{figure*}
    \centering
    \includegraphics[width=18cm, height=10cm]{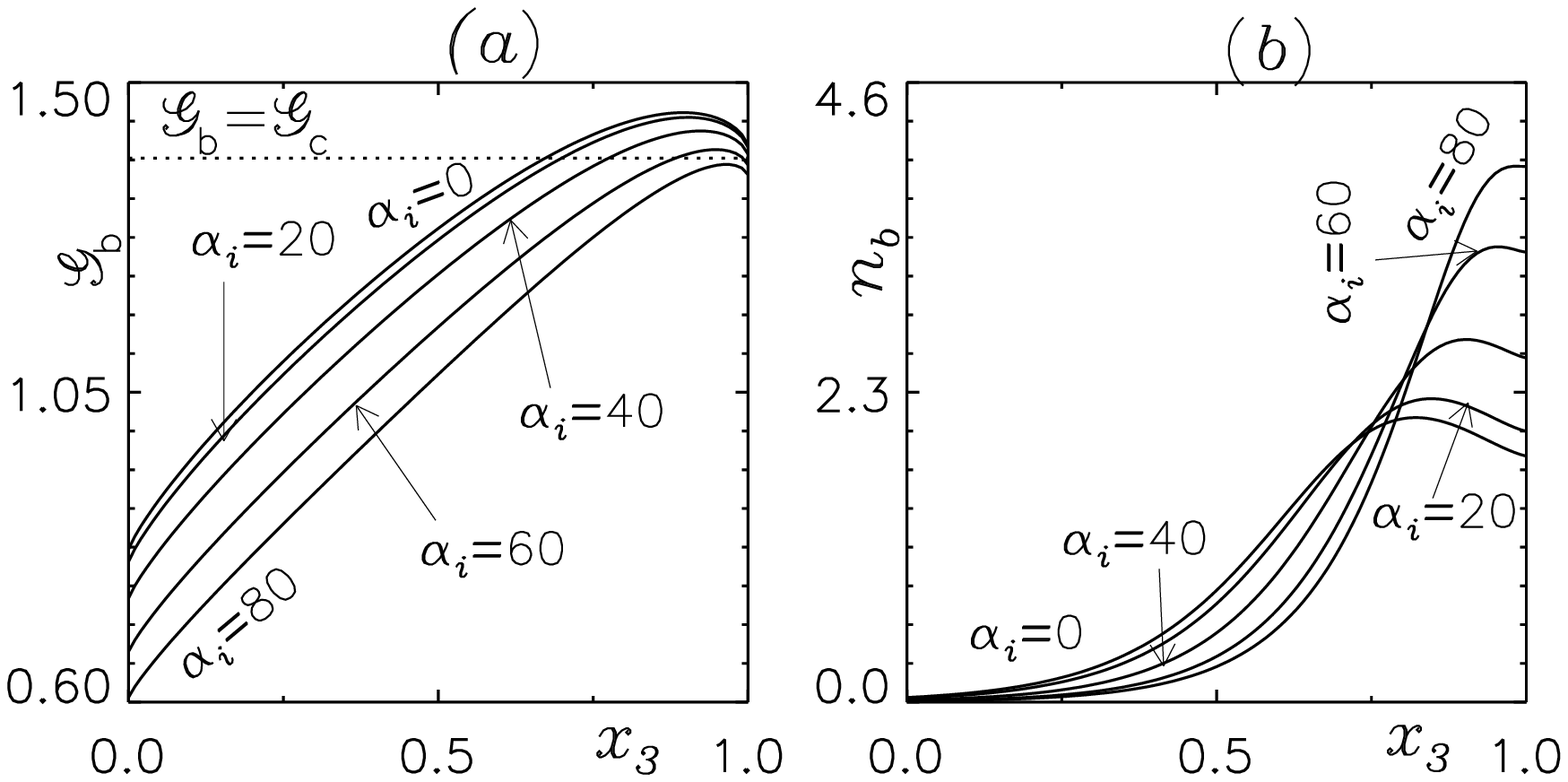}
    \caption{(a) The variation of total intensity $\mathcal{G}_b$, (b) and the corresponding basic concentration profile for $U_s=10$, $\tau_h=0.8$, $\omega=0.8$, $\mathcal{G}_c=1.39$, $A=0.2$.}
   \label{0.2A-int.eps}
 \end{figure*}
\begin{figure*}
    \centering
    \includegraphics[width=18cm, height=10cm]{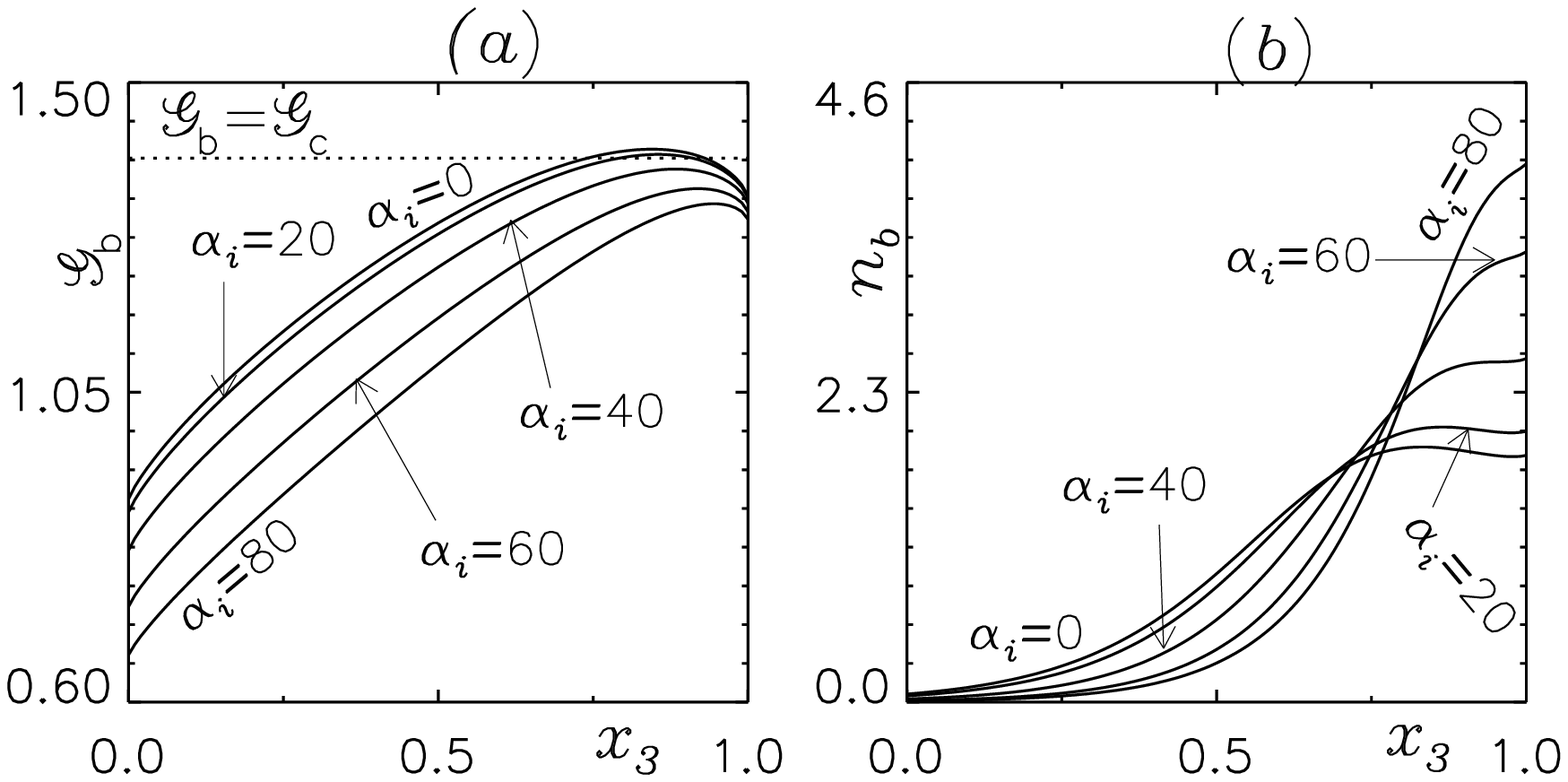}
    \caption{(a) The variation of total intensity $\mathcal{G}_b$, (b) and the corresponding basic concentration profile for $U_s=10$, $\tau_h=0.8$, $\omega=0.8$, $\mathcal{G}_c=1.39$, $A=0.78$.}
   \label{0.78A-int.eps}
 \end{figure*}
\begin{figure*}
    \centering
    \includegraphics[width=18cm, height=10cm]{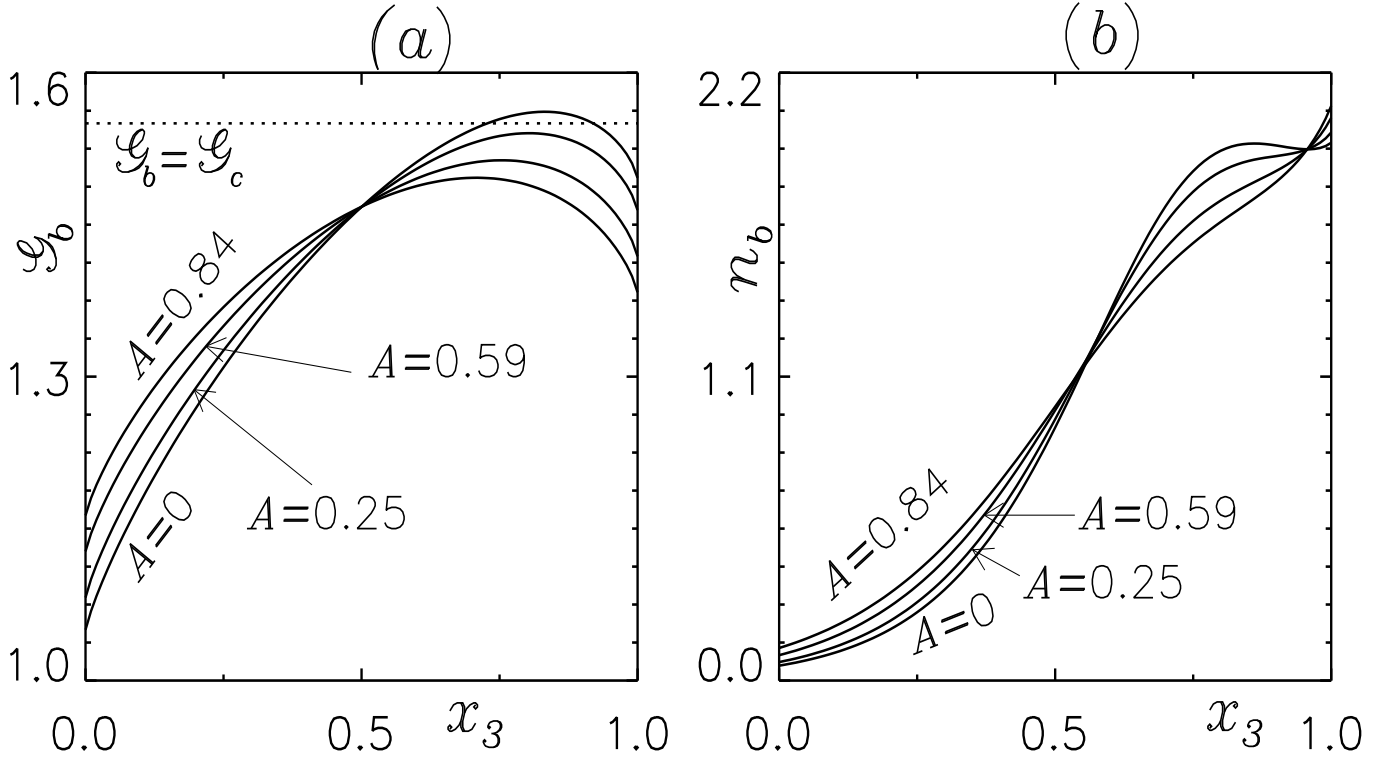}
    \caption{(a) The variation of total intensity $\mathcal{G}_b$, (b) and the corresponding basic concentration profile for $U_s=10$, $\tau_h=0.5$, $\mathcal{G}_c=1.55$, $\omega=1.0$, $\alpha_i=10$.}
   \label{10v0.5k1.55o-int.eps}
 \end{figure*}
\begin{figure*}
    \centering
    \includegraphics[width=18cm, height=10cm]{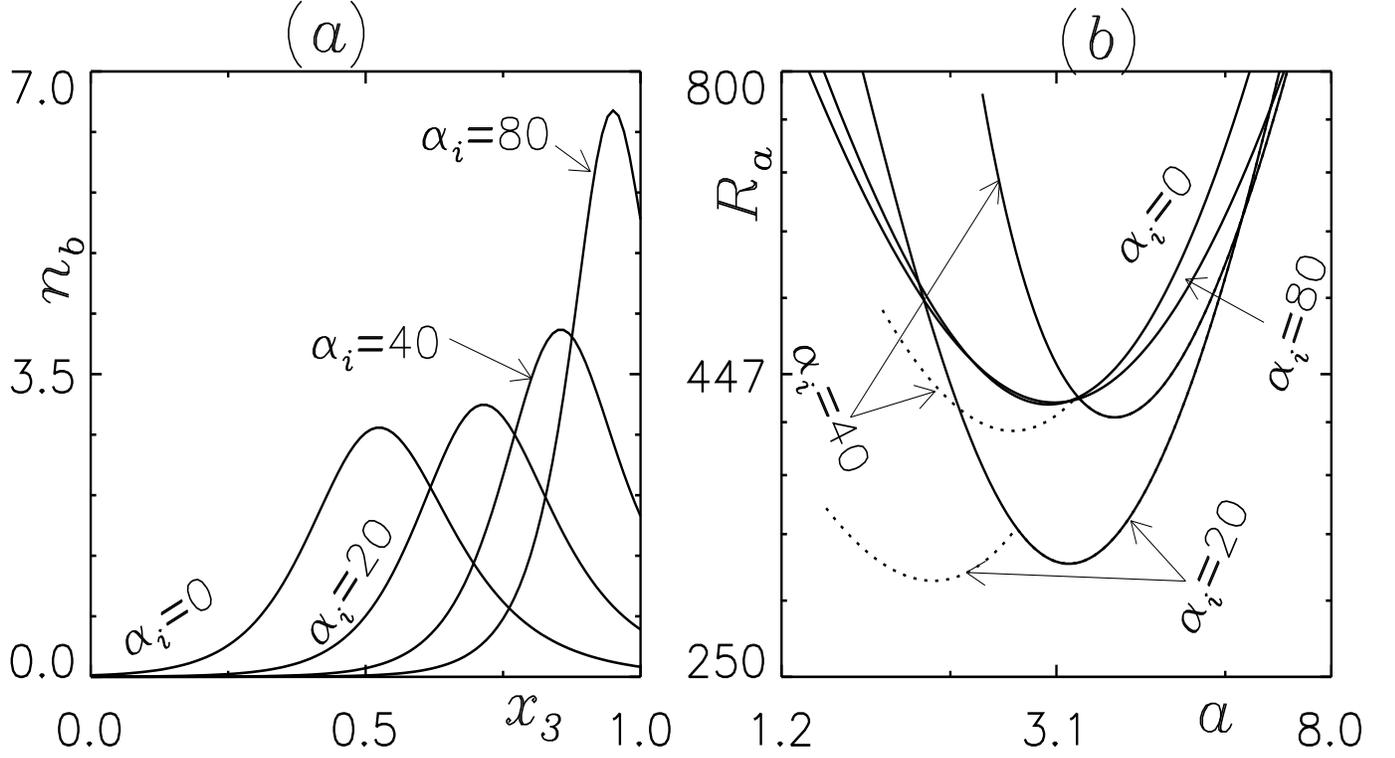}
    \caption{(a) Effect of oblique irradiation on the basic concentration and (b) the corresponding neutral curve  for $U_s=20$, $\tau_h=1.0$, $\mathcal{G}_c=1.0$, $\omega=0.605$, $A=0.38$.  In graph (b), the solid lines represent the stationary branch and the dotted lines represent the oscillatory branch.}
   \label{20v1.0k.eps}
 \end{figure*}
\begin{figure*}
    \centering
    \includegraphics[width=18cm, height=20cm]{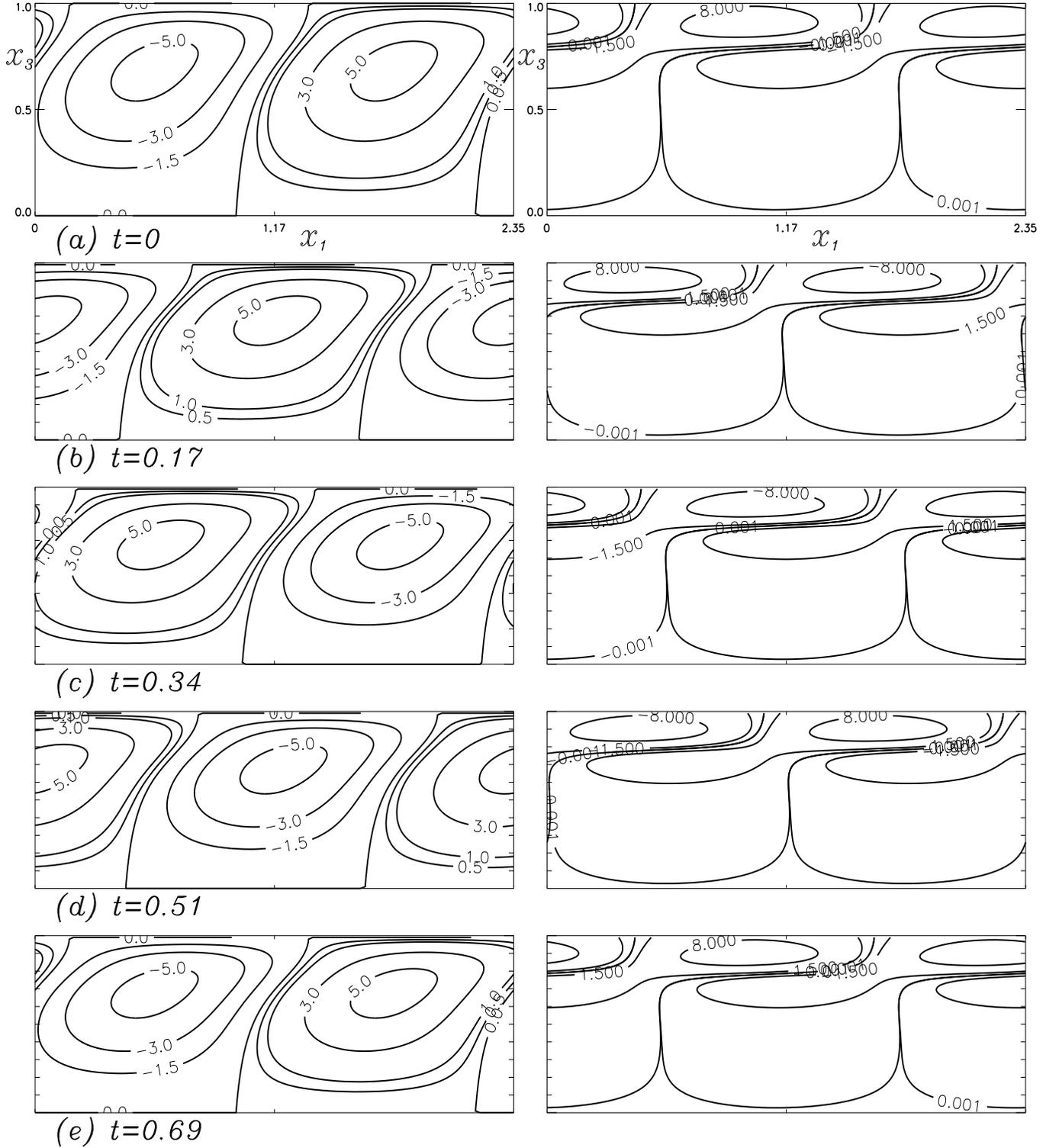}
    \caption{Flow pattern produced by the perturbed velocity $w^*$ and cell concentration $n^*$ throughout one cycle of oscillation for $U_s=20$, $\tau_h=1.0$, $\mathcal{G}_c=1.0$, $\omega=0.605$, $\alpha_i=40$, $A=0.38$, $R_a^c=401.01$, $a^c=2.67$, $Im(\gamma)=9.08$ (a) time $t=0$, (b) $t=0.17$, (c) $t=0.34$, (d) $t=0.51$, (e) $t=0.69$.}
   \label{periodic.eps}
 \end{figure*}
\begin{figure*}
    \centering
    \includegraphics[width=18cm, height=10cm]{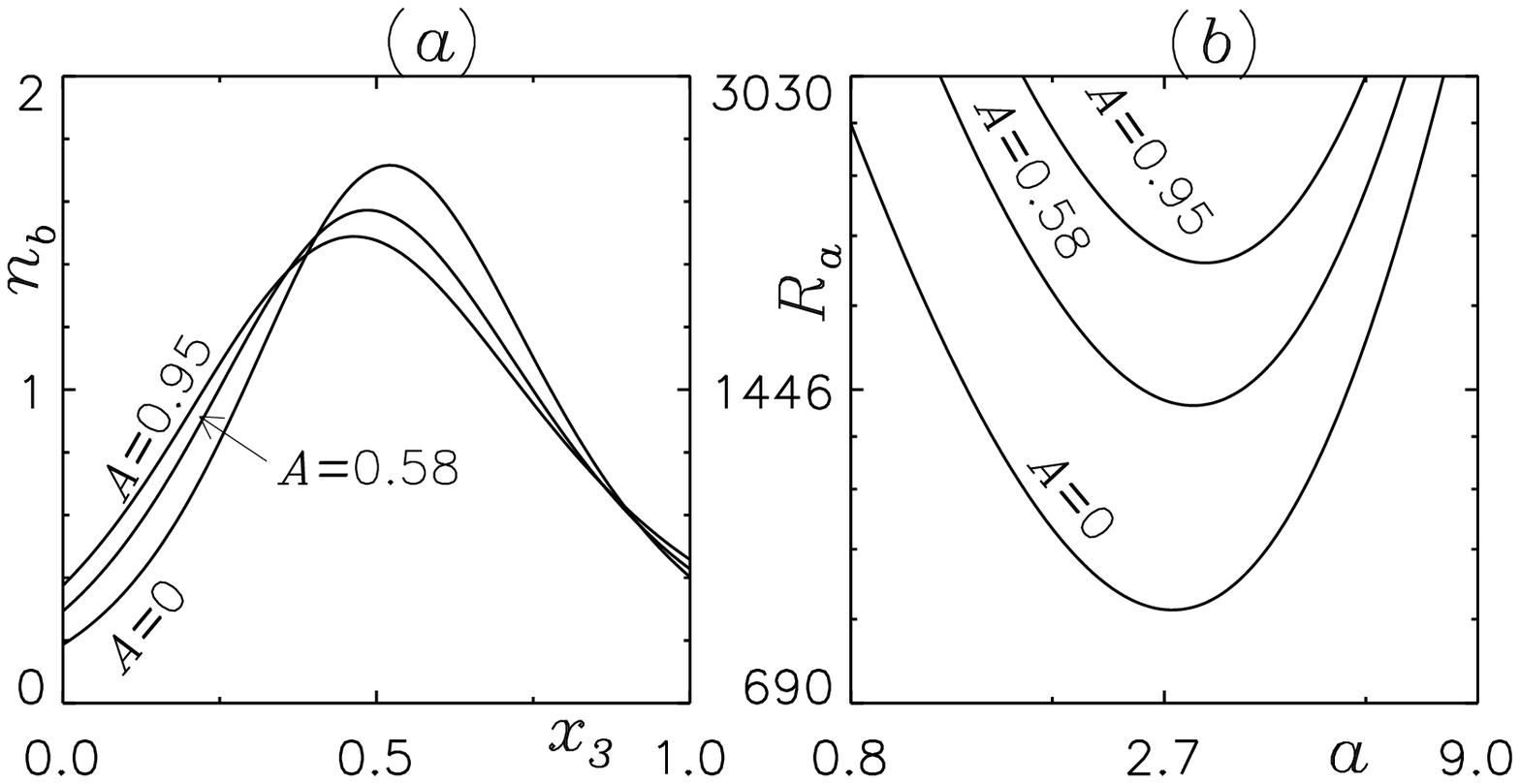}
    \caption{(a) Effect of forward anisotropic scattering coefficient on the basic concentration and (b) the corresponding neutral curve  for $U_s=16$, $\tau_h=0.5$, $\omega=0.475$, $\mathcal{G}_c=1.0$, $\alpha_i=0$.}
   \label{0.5k160th.eps}
 \end{figure*}
\begin{figure*}
    \centering
    \includegraphics[width=18cm, height=10cm]{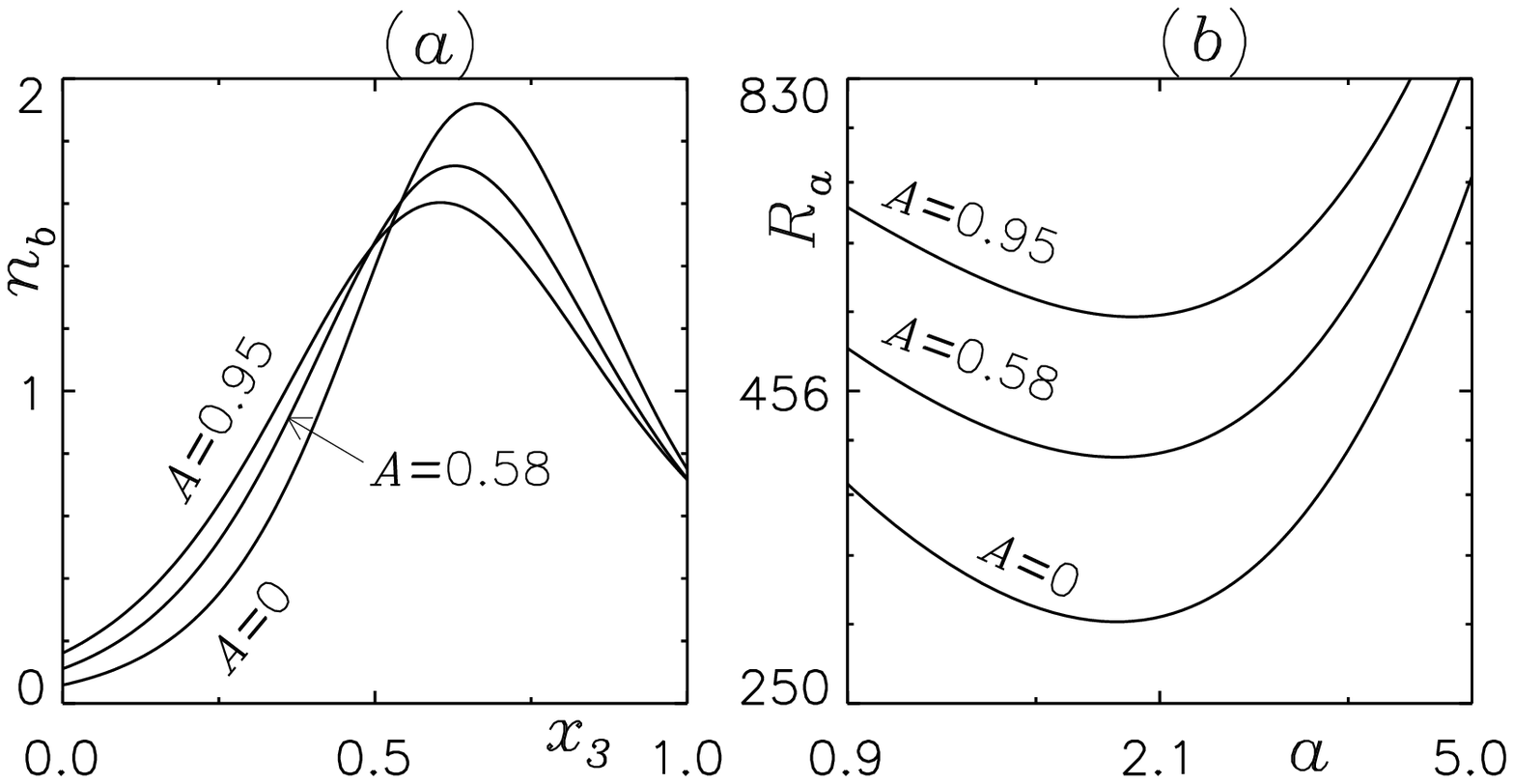}
    \caption{(a) Effect of forward anisotropic scattering coefficient on the basic concentration and (b) the corresponding neutral curve  for $U_s=16$, $\tau_h=0.5$, $\omega=0.475$, $\mathcal{G}_c=1.0$, $\alpha_i=30$.}
   \label{0.5k16v30th.eps}
 \end{figure*}
\begin{figure*}
    \centering
    \includegraphics[width=18cm, height=10cm]{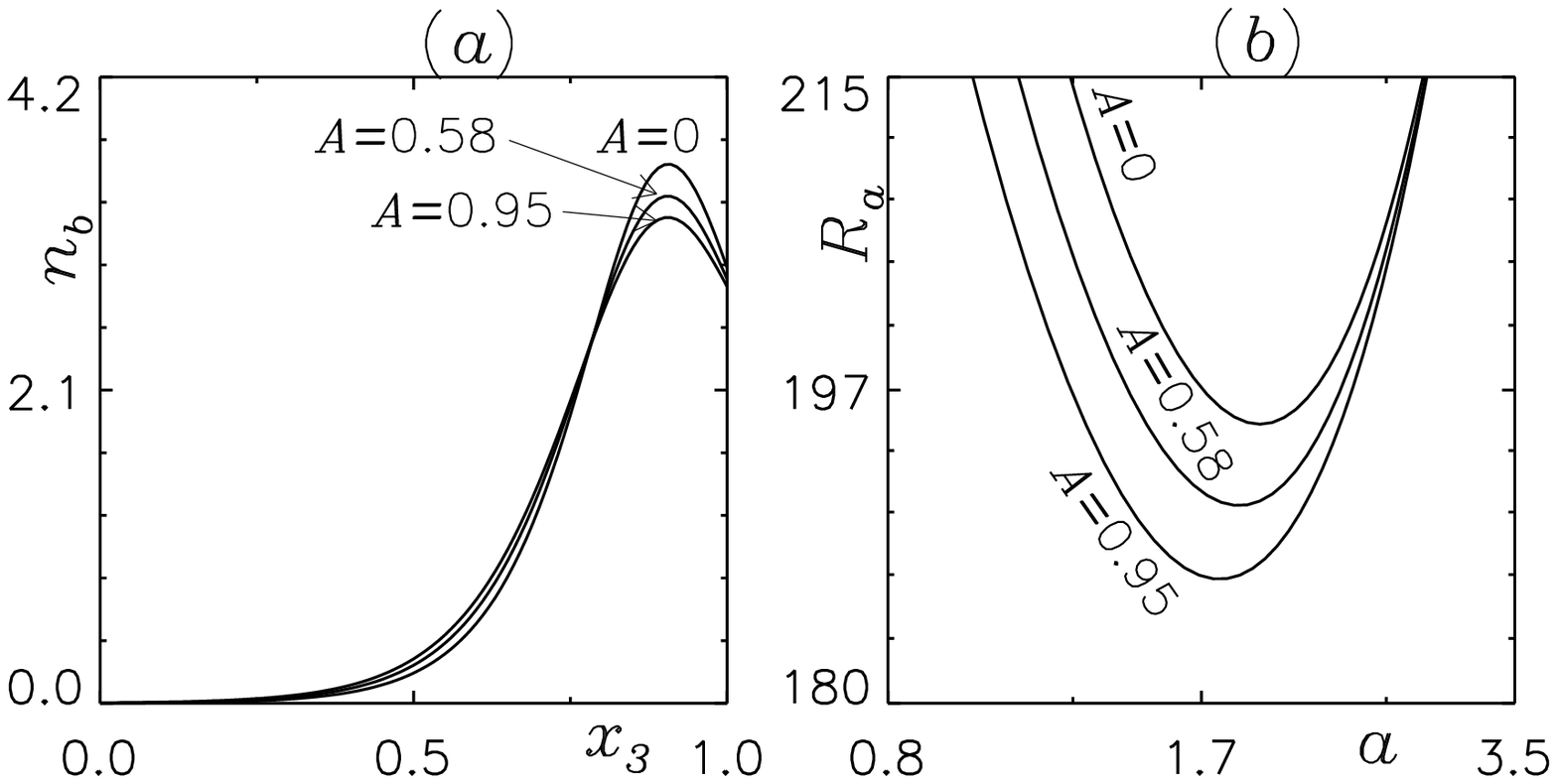}
    \caption{(a) Effect of forward anisotropic scattering coefficient on the basic concentration and (b) the corresponding neutral curve  for $U_s=16$, $\tau_h=0.5$, $\omega=0.475$, $\mathcal{G}_c=1.0$, $\alpha_i=80$.}
   \label{0.5k16v80th.eps}
 \end{figure*}
\begin{figure*}
    \centering
    \includegraphics[width=18cm, height=10cm]{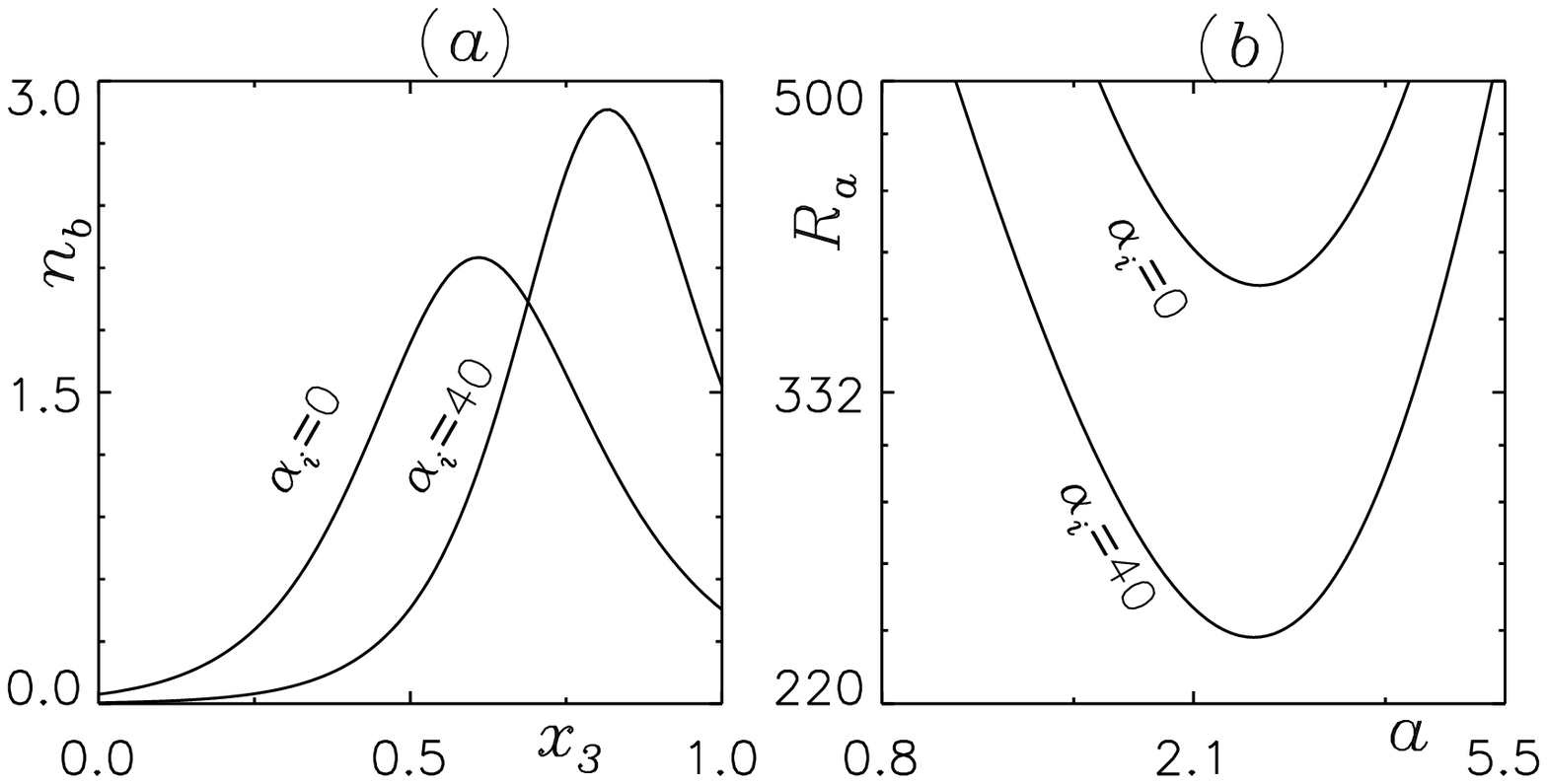}
    \caption{(a) Effect of the critical total intensity on the basic concentration and (b) the corresponding neutral curve  for $U_s=13$, $\tau_h=1.0$, $\omega=0.59$, $A=0.2$, $\mathcal{G}_c=1.0$.}
   \label{G=1.0.eps}
 \end{figure*}
\begin{figure*}
    \centering
    \includegraphics[width=18cm, height=10cm]{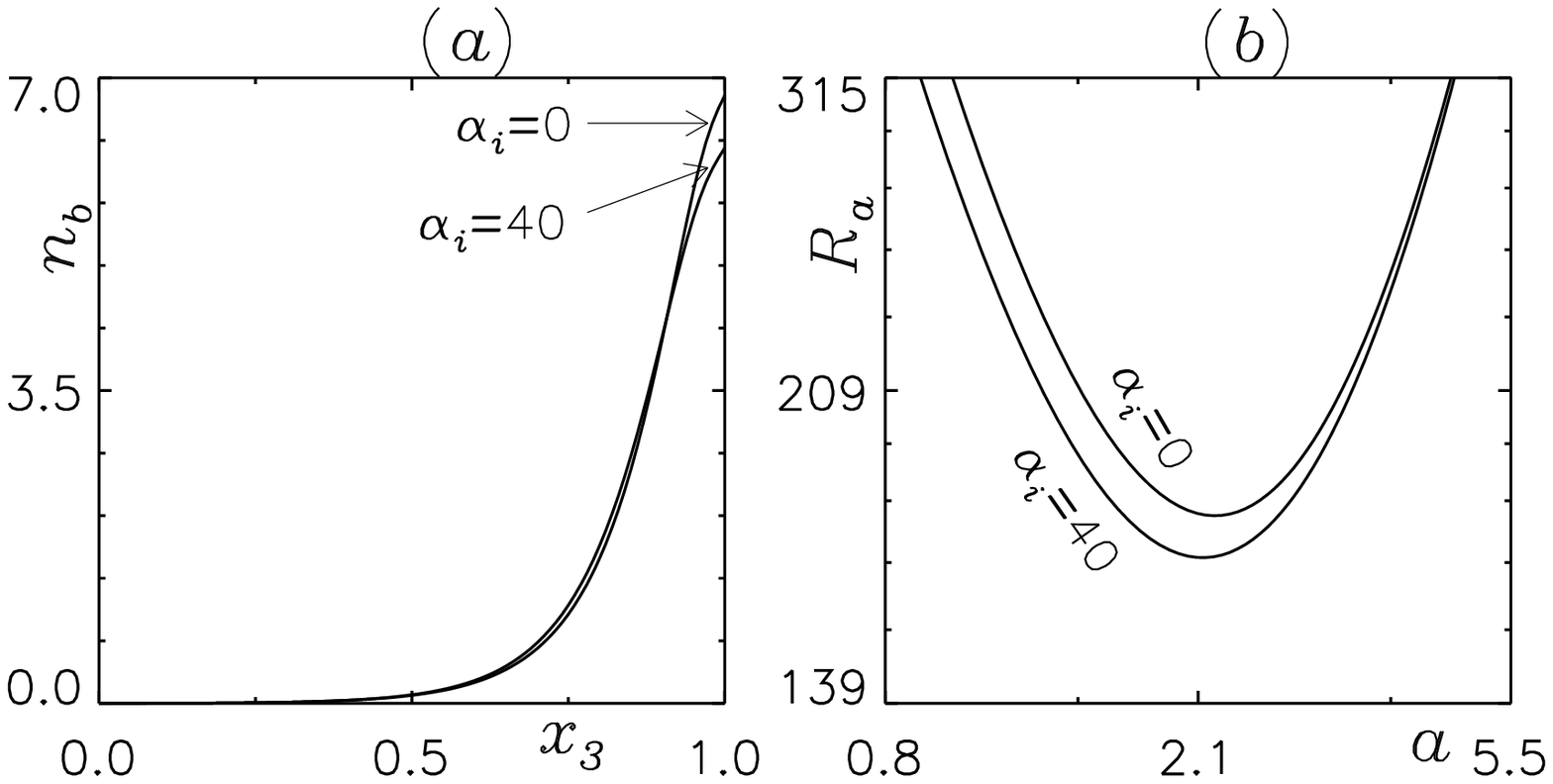}
    \caption{(a) Effect of the critical total intensity on the basic concentration and (b) the corresponding neutral curve  for $U_s=13$, $\tau_h=1.0$, $\omega=0.59$, $A=0.2$, $\mathcal{G}_c=1.39$.}
   \label{G=1.39.eps}
 \end{figure*}
\begin{figure*}
    \centering
    \includegraphics[width=18cm, height=10cm]{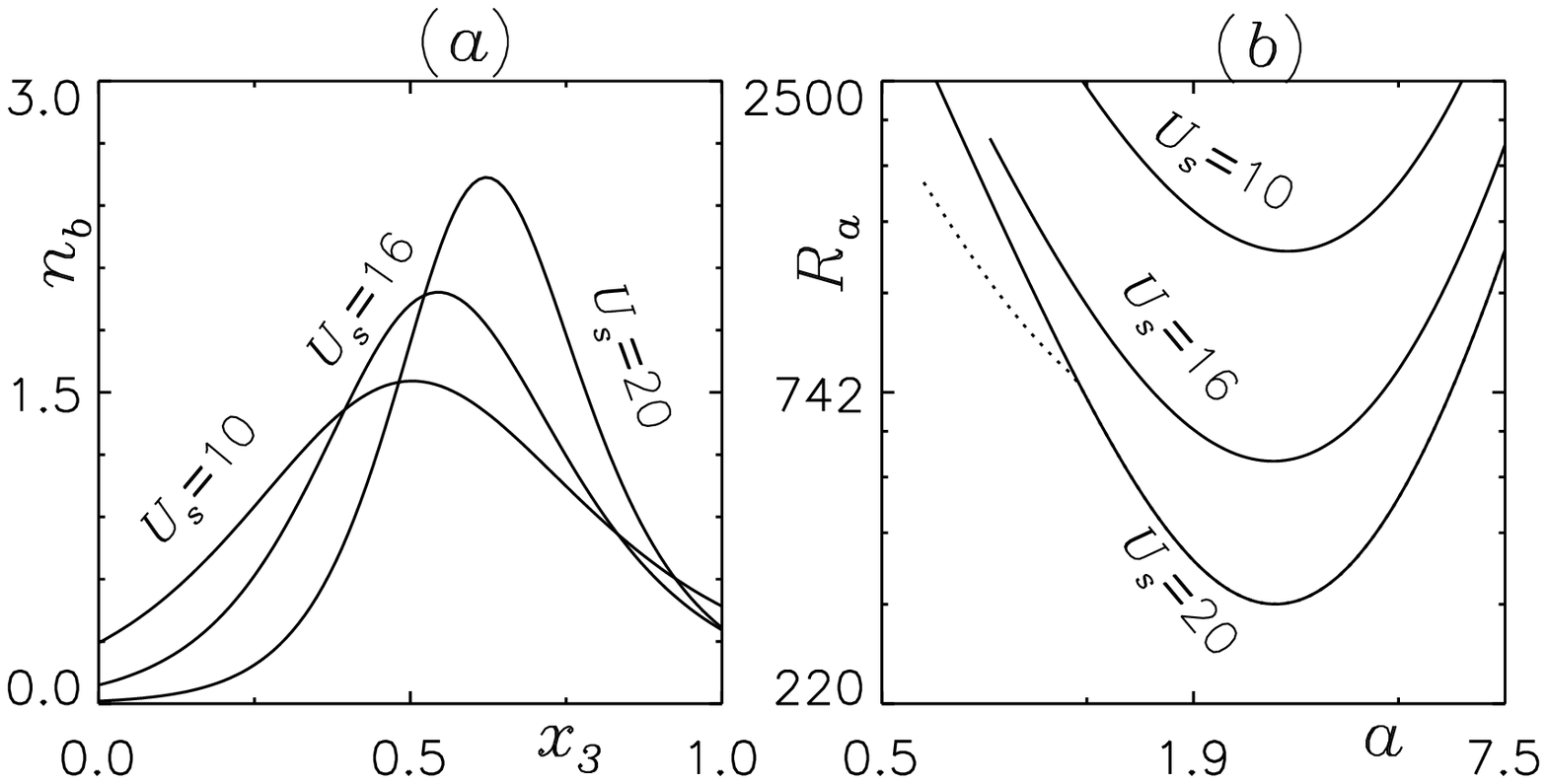}
    \caption{(a) Effect of the cell swimming speed on the basic concentration and (b) the corresponding neutral curve  for $\tau_h=0.79$, $\omega=0.55$, $A=0.38$, $\alpha_i=0$, $\mathcal{G}_c=1.0$, $\alpha_i=0$.}
   \label{th=0.eps}
 \end{figure*}
\begin{figure*}
    \centering
    \includegraphics[width=18cm, height=10cm]{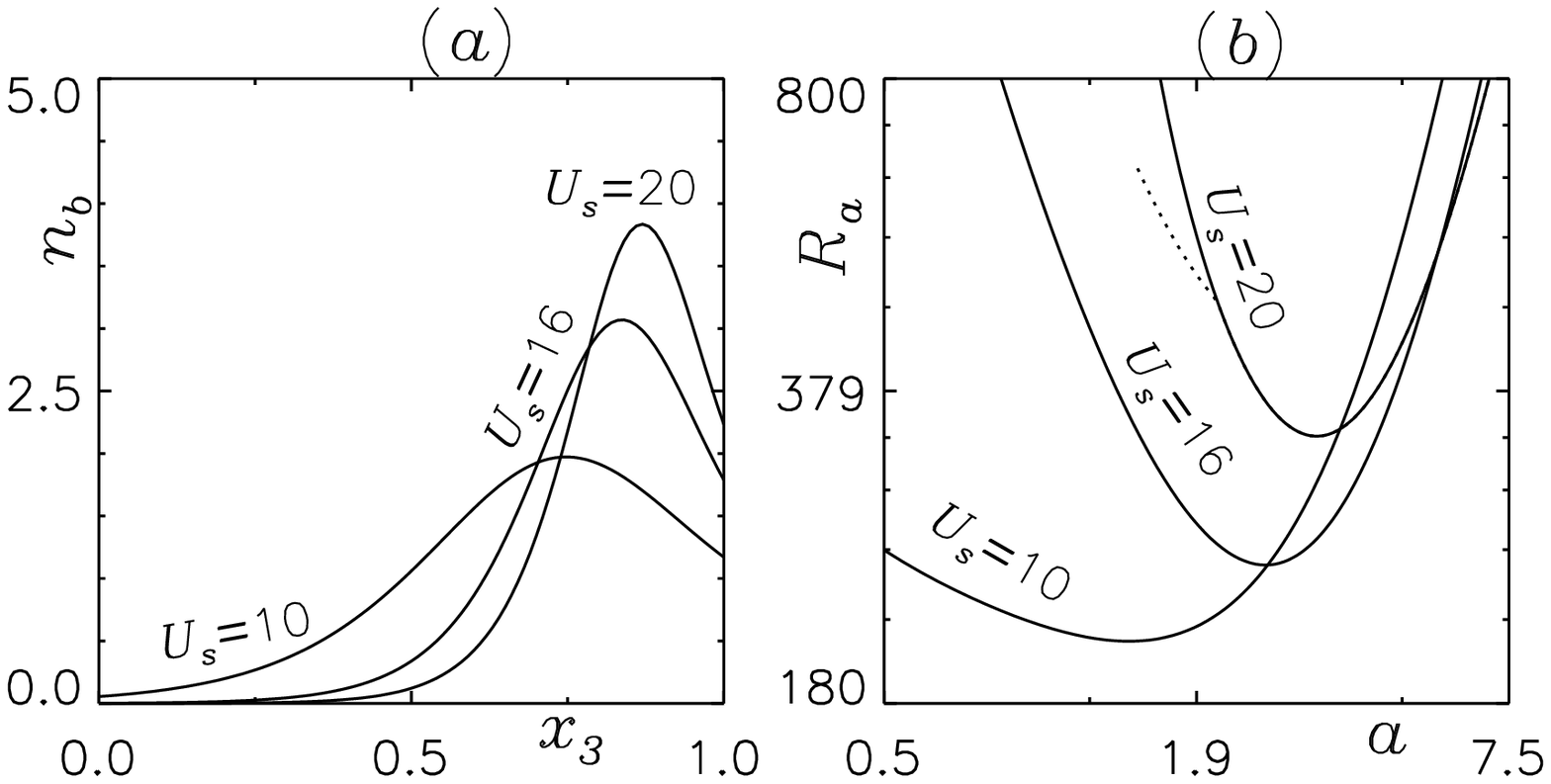}
    \caption{(a) Effect of the cell swimming speed on the basic concentration and (b) the corresponding neutral curve  for $\tau_h=0.79$, $\omega=0.55$, $A=0.38$, $\alpha_i=0$, $\mathcal{G}_c=1.0$, $\alpha_i=45$.}
   \label{th=45.eps}
 \end{figure*}

\begin{table}
\caption{Typical parameters for the suspension of \textit{Chlamydomonas} phototactic microorganism. \cite{ref9,ref13,ref17}}
\begin{tabular}{ p{5cm} p{3cm}}
\hline
\hline
Scaled average swimming speed &$U_s=20H$\\
Kinematic viscosity &$\nu=10^{-2}$cm$^2$/s\\
Schmidt number &$S_c=20$ \\
Average concentration &$\bar{n}=10^6$cm$^{-3}$\\
Average cell swimming speed &$U_c =10^{-2}$cm/s\\
Cell volume &$\vartheta=5\times10^{-10}$cm$^3$\\
Cell diffusivity &$D=5\times10^{-4}$cm$^2$/s\\
Ratio of cell density &$\Delta \varrho/\varrho=5\times 10^{-2}$ \\
Cell radius &$a=10^{-3}$cm\\
\hline
\hline
\end{tabular}\\

\label{tab:table1}
\end{table}

\section{Numerical solution}

The system of ordinary differential equations (\ref{eqn:71})-(\ref{eqn:75}) is of the seventh order with seven boundary conditions. To solve an ordinary differential equation of order seven, we use the Newton-Raphson-Kantorovich iteration-based method. \citep{ref19} We investigate the linear stability of the basic state by drawing neutral curves in the $(a, R_a)$ plane. In addition, the numerical method is tested with a number of different mesh sizes and characteristics. It proves that the solutions that are created by utilizing the technique for the same parameters on multiple meshes are and always agree to five significant figures or more for a minimum of $51$ mesh points and that the solutions are independent of both time and grid. Points where the real component of the growth rate $Re(\gamma)$ is zero constitute a neutral curve. The existence of an overstable or oscillatory solution is possible if the imaginary component of the growth rate $Im(\gamma)$ is not zero on the neutral curve. On the other hand, if the imaginary part of the growth rate $Im(\gamma)$ is zero on such a curve, then the perturbation to the basic state is stationary, and the principle of exchange of stabilities holds. \citep{ref12} The growth rate $Re(\gamma)$ and a hypothetical oscillation with frequency $Im(s)/2\pi$ are both accounted for in the time dependency $\exp(\gamma t) = exp[(Re(\gamma)+i Im(\gamma))t]$. For stable layering, all wave vectors have a negative growth rate $Re(\gamma)$. For a narrow range of wave numbers, just above the convective beginning, $Re(\gamma)$ becomes positive. The neutral curve $R_a^{n}(a) (n = 1, 2, 3,... )$ has an unlimited number of branches for a given set of other parameter ranges, each of which represents a unique solution to the linear stability problem. The one on which $R_a$ has its smallest value, $R_a^c$, is the most interesting branch of the possible solutions. It has been shown that the pair $(a^c, R_a^c)$ represents the most unstable solution. With this information, one may calculate the wavelength of the original perturbation using the formula $\lambda^c = 2\pi/a^c$. The bioconvective solutions are made up of individual convection cells that are layered one on top of the other throughout the length of the suspension. If a solution has $n$ convection cells that are piled vertically one on top of the other, then it is said to be of mode $n$. The most unstable solution, which happens to be mode $1$, may frequently be found on the $R^1(a)$ branch of the neutral curve. This is the case in many different situations. 

We presume that we are working with a phototactic microorganism comparable to \textit{Chlamydomonas} to determine the parameters needed for the current study. To make the model more reasonable and similar to other research on phototactic bioconvection, we adopt the same parameter values as References \citet{ref9,ref15,ref16a,ref17} (see Table \ref{tab:table1}). The radiation parameters required here are computed the same as in \citet{ref15}. Therefore, the optical depth ranges from $0.25$ to $1$ for a $0.5$ cm depth suspension. For a suspension depth of $0.5$ cm, $U_s = 10$ is the corresponding scaled swimming speed, and for a depth of $1.0$ cm, $U_s = 20$ is the corresponding scaled swimming speed (see Table \ref{tab:table1}). The range of the incident angle is the same as $0^{\circ}\le\alpha_0\le 80^{\circ}$ given in ~\citet{ref17,ref18}. To facilitate a comparison of our model with other rational representations of phototactic bioconvection, the parameters $S_c = 20$ and $I^0 = 1.0$ have been held constant throughout.

The taxis function calculates the critical total intensity $\mathcal{G}_c=1.39$ for $\Upsilon=0.4$. Figure (\ref{0.2A-int.eps}) illustrates the effect of oblique irradiation on total intensity and corresponding basic concentration for $U_s=10$, $\tau_h=0.8$, $\omega=0.8$, $\mathcal{G}_b=1.39$, $A=0.2$. For $\alpha_i=0$, the location of critical total intensity is around $x_3=0.67$, and algae cells accumulate near $x_3=0.81$. The cells that are located above the location of $\mathcal{G}_c$ have a phototactic response that is negative, whereas the cells that are located below the location of $\mathcal{G}_c$ have a phototactic response that is positive. Due to the self-shading becoming dominant at a large angle of incidence, microorganisms acquire low intensities at a constant interior suspension depth. Thus, the position of the maximal basic concentration shifts toward the top of the domain as the value of $\alpha_i$ grows up to $80$. Figure (\ref{0.78A-int.eps}) illustrates the effect of oblique irradiation on total intensity and corresponding basic concentration for $U_s=10$, $\tau_h=0.8$, $\omega=0.8$, $\mathcal{G}_c=1.39$, $A=0.78$. For $\alpha_i=0$, $\mathcal{G}_c$ may be found in a uniform suspension at two different depths: $x_3=0.73$ and $x_3=0.92$. Hence, cells that are above $x_3=0.92$ and below $x_3=0.73$ are considered to be positively phototactic, whereas cells that are in the middle are considered to be negatively phototactic. The scattering that occurs in the suspension is the cause of the accumulation of cells in two distinct locations. Therefore, the cells accumulate around $x_3 = 0.82$ and the top of the suspension. The total intensity throughout the whole domain drops below the critical intensity as the angle of incidence increases. This causes the entire suspension to become positively phototactic, and the cells to gather in a single location. As a consequence, the location of the maximal basic concentration is condensed into a single point and moved to the top of the suspension. 

The variation in total intensity, $\mathcal{G}_b$, and basic state concentration that occurs over the depth of a uniform suspension is depicted in Figure (\ref{10v0.5k1.55o-int.eps}), with $U_s=10$, $\tau_h=0.5$, $\mathcal{G}_c=1.55$, $\omega=1.0$, $\alpha_i=10$. In this case, the upper half of the suspension experiences a drop in $\mathcal{G}_b$ with increasing $A$, whereas the bottom half has the reverse trend. This is due to the fact that a larger value of A causes the forward scattering to scatter a greater amount of energy into the forward directions. In addition to this, the change in $\mathcal{G}_b$ during the suspension does not follow a monotonic pattern. As a result, the critical total intensity might happen in two distinct locations inside the suspension. For $\alpha_i=0$, $\mathcal{G}_c$ may be found in a uniform suspension at two different depths: $x_3=0.72$ and $x_3=0.9$. Thus, the cells accumulate around $x_3 = 0.85$ as well as the top of the suspension. As the anisotropic scattering coefficient $A$ increases, the critical intensity shifts to a single location, and the cells only exhibit positive phototaxis (see Figure \ref{10v0.5k1.55o-int.eps}(b)). 

In Figure (\ref{20v1.0k.eps}), the base concentration profiles and the related neutral curves on varying the angle of incidence are illustrated as $\alpha_i = 0, 20, 40, 60$, and $80$, respectively. This is done while keeping the other governing parameters, $U_s=20$, $\tau_h=1.0$, $\mathcal{G}_c=1.0$, $\omega=0.605$, $A=0.38$, constant. If $\alpha_i=0$, the highest concentration of the basic substance occurs in the middle of the domain. As the most unstable solution is in the stable branch of the corresponding neutral curve, the disturbance in the initial state is stable for $\alpha_i=0$. The bioconvective solution is overstable if $\alpha$ is raised to $20$ due to the most unstable mode persisting on the oscillatory branch of the neutral curve. The patterned behavior while $\alpha_i=40$ is similar to the case of $\alpha_i=20$. Furthermore, in this instance, a single oscillatory branch bifurcates from the corresponding stationary branch of the neutral curve at around $a=3.34$, however, the oscillatory branch retains the most unstable bioconvective solution (see Figure \ref{20v1.0k.eps}(b)). As a result, overstability occurs at $a^c=2.67$ and $R_a^c=401.01$. $\gamma=0\pm 9.08 i$ are the two eigenvalues that are determined to be complex conjugates at this point. The change that can be seen here is referred to as a Hopf bifurcation. It can be shown that the bioconvective flow patterns that correspond to the complex conjugate pair of eigenvalues are mirror images of one another. The oscillation period is $2\pi/Im(\gamma)=0.69$ units. On a timeframe that is considerably shorter than the projected period of overstability, the bioconvective fluid movements transition into a fully nonlinear state. Hence, the perturbed eigenmodes $w^*$ and $n^*$ may be used to view the convection cells and flow patterns throughout one cycle of oscillation (see Figure \ref{periodic.eps}). It demonstrates that a traveling wave solution is progressing in the direction of the figure's left side. 

The impacts of the anisotropic scattering coefficient, $A$, on the basic concentration and neutral curves for a variety of incidence angles are depicted in Figures \ref{0.5k160th.eps}–\ref{0.5k16v80th.eps} for the fixed parameters $U_s=16$, $\tau_h=0.5$, $\omega=0.475$, $\mathcal{G}_c=1.0$. At $A=0$, the maximum concentration of the basic state for $\alpha_i=0$ is located at the domain's mid-height (see Figure \ref{0.5k160th.eps}). As $A$ is raised to $0.58$ and $0.95$, respectively, the position of the greatest base concentration moves away from the top of the suspension, and the thickness of the upper stable layer rises in a consistent manner. Hence, the influence of buoyancy, which has the tendency to restrict convective fluid motion, grows monotonically as A is raised, and as a consequence, the critical wavenumber and the critical Rayleigh number increase as well. Figure \ref{0.5k16v30th.eps} shows the basic concentration and the corresponding neutral curves for $\alpha_i=30$. For $A=0$, the position of the maximum basic state concentration occurs at $x_3=0.66$. In this case, the maximum basic concentration decreases as the anisotropic scattering coefficient $A$ increases, and their position moves toward the bottom of the suspension. As a result, the critical wavenumber and the critical Rayleigh number increase. The maximum concentration at the basic state at $\alpha_i=30$ for all values of $A$ is comparatively higher than the case of $\alpha_i=0$ and also the thickness of the upper stable layer is lower. As a result, The critical Rayleigh and wavenumber decrease for $\alpha_i=30$ compared to $\alpha_i=0$ for all values of $A$. Figure \ref{0.5k16v80th.eps} shows the basic concentration and the corresponding neutral curves for $\alpha_i=80$. In this case, the maximum basic concentration decreases as the anisotropic scattering coefficient $A$ increases. As the value of $A$ is varied from higher to lower, the concentration gradient in the top stable region becomes steeper. As a result, the critical Rayleigh number decreases as $A$ is increased from $0$ to $0.95$. The impacts of the incidence angles, $\alpha_i$, on the basic concentration and neutral curves for a variety of the critical total intensity $\mathcal{G}_c$ are depicted in Figures \ref{G=1.0.eps}–\ref{G=1.39.eps} for the fixed parameters $U_s=13$, $\tau_h=1.0$, $\omega=0.59$, $A=0.2$. Figure \ref{G=1.0.eps} shows the basic concentration and the corresponding neutral curves for $\mathcal{G}_c=1.0$. In this case, the maximum basic concentration occurs around $x_3=0.6$ for $\alpha=0$. As $\alpha_i$ increases to $40$, the maximum basic concentration rises and their position shifts towards the top of the suspension. As a result, the critical Rayleigh number decreases as $\alpha_i$ is increased from $0$ to $40$. Figure \ref{G=1.39.eps} shows the basic concentration and the corresponding neutral curves for $\mathcal{G}_c=1.39$. The maximum concentration at the basic state at $\mathcal{G}_c=1.39$ for all values of $\alpha_i$ is comparatively higher than the case of $\mathcal{G}_c=1.0$ and also the thickness of the lower unstable layer is higher. Thus, for the higher critical total intensity, the critical Rayleigh number and wavenumber decrease, and the system becomes more unstable.

The impacts of the cell swimming speed, $U_s$, on the basic concentration and neutral curves for a variety of the incidence angle are depicted in Figures \ref{th=0.eps}–\ref{th=45.eps}. Figure \ref{th=0.eps} shows the basic state concentration and the neutral curves for the fixed parameters $\tau_h=0.79$, $\omega=0.55$, $A=0.38$, $\alpha_i=0$, $\mathcal{G}_c=1.0$. At $U_s = 10$, the highest base concentration occurs near the middle of the suspension. As the cell swimming speed $U_s$ increases to $16$, $20$, the maximum basic concentration rises, and their position shifts towards the top of the suspension. The width of the upper stable layer monotonically decreases as $U_s$ is varied as $10$, $16$, and $20$ respectively. As a result, the critical wavenumber and the critical Rayleigh number decrease for higher cell swimming speed. For $U_s=10$ and $U_s=16$, the most unstable solutions are on the stationary branch. For $U_s=20$, a single oscillatory branch bifurcates from the stationary branch at wavenumber $a= 1.18$. However, the most unstable solution still remains on the stationary branch. For $\alpha_i=45$, the basic state concentration and the neutral curves are depicted in Figure \ref{th=45.eps}. Here, $x_3=0.75$, $x_3=0.83$, and $x_3=0.87$ are the locations of the greatest base concentration for $U_s=10,16$, and $20$, respectively. A steep concentration in the basic steady state implies a higher concentration gradient which supports the bioconvection. But at a higher swimming speed, the positive phototaxis offers higher resistance to the cells residing in the bioconvective plume. Hence, the latter effect dominates the higher gradient leading to the higher critical Rayleigh number for a higher swimming speed.


\section{Conclusion}
The effect of oblique collimated irradiation on the suspension of phototactic bioconvection with forward anisotropic scattering is explored in this work. We considered the bottom layer of the suspension to be rigid while the top layer to be stress-free. This model has been used to investigate the linear stability of the suspension. It is observed that self-shadowing in the basic state dominates through an increase in the slant-path length upon an increase in the incident angle. Therefore, due to variations in the angle of incidence at a constant internal depth, the algae cells receive light of low intensity. As a consequence of this, the position of the maximum basic concentration moves closer to the top of the suspension if there is an increase in the angle of incidence, and the value of the maximum basic concentration rises accordingly. At larger values of the forward scattering coefficient $A$, there is a decrease in the total intensity in the upper region of the suspension. On the other hand, the converse is true for the lower region. This is due to the fact that larger values of the forward scattering coefficient $A$ cause the forward scattering to scatter a greater amount of energy in the forward direction. As a consequence of this, the position of the critical total intensity, which is found in the top region of a uniform suspension, is higher or equal for a greater value of the forward scattering coefficient $A$ than it is for a value of $A$ that is lower. At the lower part of the suspension, one would find the opposite to be true. Moreover, as the value of A increases, the variance of the total intensity becomes less steep across the suspension. The critical Rayleigh number is influenced by a number of factors, including the thickness of the upper stable layer, the maximum concentration of the basic state, and the concentration gradient in the upper stable layer, among other things. Phototaxis either supports or prevents the convection in the suspension. Also, The gravitationally unstable domain below the sublayer supports the convection whereas the domain above prevents it. Therefore, Oscillatory solutions are observed as a result of the competition between stabilizing and destabilizing processes. \cite{ref10} For higher critical total intensity, the location of the maximum basic concentration is at the top of the suspension. In this case, the marginal state is stationary and the critical Rayleigh number is lower for a higher value of the incidence angle. As the value of the critical total intensity decreases, the location of the maximum basic concentration shifts away from the top of the suspension. Thus, the critical Rayleigh and wavenumber increase as compared to the case of higher critical total intensity, and the system becomes more stable. The cells’ swimming speed affects the suspension; at faster speeds, oscillatory or overstable solutions can be found. In the proposed phototaxis model, theoretical predictions should be compared to quantified experimental findings on bioconvection in a phototactic algae solution. To our regret, there are currently no statistics of this kind available. This
is due to the fact that it would be necessary to find a type of microbe that is mostly phototactic, while most species of algae in the natural environment are also gravitactic or gyrotactic.\cite{ref8} Further experimentation is required to calculate optical depths, phototaxis functions, and diffusion coefficients in order to learn more about these intriguing phenomena.


\section{Acknowledgements}
This study was supported by the University Grants Commission, Grants number 191620003662, New Delhi (India). 
\nocite{*}
\bibliography{aipsamp}


\end{document}